\documentclass[10pt]{amsart}

\usepackage{amsmath}
\usepackage{amssymb}
\usepackage{amsfonts}
\usepackage{graphicx}
\usepackage{epsfig}

\newtheorem{proposition}{Proposition}[section]
\newtheorem{theorem}{Theorem}
\newtheorem{lemma}[proposition]{Lemma}
\newtheorem{remark}[proposition]{Remark}
\newtheorem{definition}[proposition]{Definition}  
\newtheorem{corollary}[proposition]{Corollary}  
\newtheorem{notation}[proposition]{Notation}

\def\C{{\mathbb C}}

\def\R{{\mathbb R}}
\def\Z{{\mathbb Z}}
\def\N{{\mathbb N}}
\def\S1{{\mathbb S}^1}

\def\T{{\R / \Z}}
\def\cyl{\T \times \R}
\def\Cstar{\C^{*}}

\def\cA{{\mathcal A}}
\def\cC{{\mathcal C}}
\def\cD{{\mathcal D}}

\def\cJ{{\mathcal J}}
\def\cL{{\mathcal L}}
\def\cM{{\mathcal M}}
\def\cO{{\mathcal O}}

\def\cS{{\mathcal S}}
\def\cT{{\mathcal T}}

\def\cU{{\mathcal U}}

\def\btheta{{\bar{\theta}}}
\def\bV{\bar{V}}

\def\tf{{\widetilde{f}}}
\def\tg{{\widetilde{g}}}

\def\tF{{\widetilde{F}}}

\def\tgamma{{\widetilde{\gamma}}}
\def\tmu{{\widetilde{\mu}}}
\def\trho{{\widetilde{\rho}}}
\def\ttau{{\widetilde{\tau}}}

\def\htau{{\widehat{\tau}}}
\def\hmu{{\widehat{\mu}}}

\def\thetanot{{\theta_0}}
\def\tnot{{t_0}}

\def\thetap{\theta^\prime}

\def\la{\lambda}

\def\Omegalt{\Omega_{\la, \tau}}
\def\lambdanot{\lambda_0}
\def\Alt{A_{\la,\tau}}
\def\tAlt{\tilde{A}_{\la,\tau}}
\def\Tdist{\operatorname{dist}_{\T}}
\def\ladist{\operatorname{dist}_\la}
\def\distL{\operatorname{dist}_{L}}
\def\lip{\operatorname{Lip}}
\def\bthetap{{\btheta^\prime}}
\def\multt{\operatorname{mult}_\tau}

\def\ml{{\bf m}_\ell}

\def\lainzoo{\la \in (0,1)}
\def\lainzoc{\la \in (0,1]}
\def\lainlazoo{\la \in [\la_0, 1)}

\newcommand{\dist}{\mbox{\rm dist}}

\newcommand{\m}{{\bf m}}
\newcommand{\tht}{(\theta, t)}

\def\proof{\par\noindent {\bf Proof.} }

\begin{document}

\title{On the topology of solenoidal attractors of the cylinder \\ \medskip Sur 
la topologie des attracteurs de type soleno\"{\i}de du  cylindre}

\thanks{This work has been partially supported by  Mecesup UCN-0202,  Mecesup PUC-UCH-0103 and IMPA}

\thanks{\dag Partially supported by Fondecyt \# 1010865}
\thanks{\ddag Partially  supported  by Fondecyt \#1020711.}
\thanks{$\ast$Partially supported by Fundaci\'on Andes, PROSUL,  and Fondecyt \#1070711}
\thanks{$\star$Partially  supported by Fondecyt \#1000047,  DGICT-UCN and Fundaci\'on Andes.}

\date{March 22, 2004}

\maketitle
\markboth{On the topology of solenoidal attractors of the cylinder}{On the topology of solenoidal attractors 
of the cylinder}

\begin{center}
{Rodrigo  Bam\'on\dag}\\
{Departamento de Matem\'aticas, Facultad de Ciencias, Universidad de Chile},\\  {Casilla 653, Santiago, Chile}\\
{\tt rbamon@uchile.cl}\\
{\hfill }\\
{Jan  Kiwi\ddag}\\
{Facultad de Matem\'aticas,
Pontificia Universidad Cat\'olica},\\
{Casilla 306, Correo 22, Santiago,
Chile.} \\
{\tt jkiwi@puc.cl}\\
{\hfill}\\
{Juan Rivera-Letelier$\ast$}\\
{Departamento de Matem\'aticas, Universidad Cat\'olica del Norte},\\
{Casilla  1280,
Antofagasta, Chile.}\\
{\tt juanrive@ucn.cl}\\
{\hfill} \\
{Richard  Urz\'ua$\star$}\\
{Departamento de Matem\'aticas, Universidad Cat\'olica del Norte},\\
{Casilla  1280,
Antofagasta, Chile.}\\
{\tt rurzua@ucn.cl}\\
\end{center}

\renewcommand{\abstractname}{Abstract}
\begin{abstract}
We study the dynamics of skew product endomorphisms acting on the cylinder~$\cyl$, of the form
$$
\tht \mapsto (\ell \theta, \la \theta + \tau (\theta)),
$$
where $ \ell \geq 2$ is an integer, $\la \in (0,1)$ and $\tau: \T \to \R$ is a continuous function.
We are interested on {\it topological} properties of the global attractor $\Omegalt$ of this map.
Given $\ell$ and a Lipschitz function $\tau$, we show that the attractor set $\Omegalt$ is homeomorphic to a closed topological annulus for all $\la$ sufficiently close to~$1$. 
Moreover, we prove that $\Omegalt$ is a Jordan curve for at most finitely many $\la \in (0,1)$. 

 These results rely on a detailed study of iterated ``cohomological'' equations of the form $\tau = \cL_{\la_1} \mu_1$,
$\mu_1 = \cL_{\la_2} \mu_2, \dots$, where $\cL_\la \mu = \mu \circ \ml - \la \mu$ and $\ml: \T \rightarrow \T$ denotes
the multiplication by $\ell$ map.
We show the following finiteness result: each Lipschitz function $\tau$ can be written in a canonical way as,
$$
\tau = \cL_{\la_1} \circ \cdots \circ \cL_{\la_m} \mu,
$$
where $m \ge 0$, $\lambda_1, \ldots, \lambda_m \in (0, 1]$ and the Lipschitz function $\mu$ satisfies $\mu \neq \cL_\la \rho$ for every continuous function $\rho$ and every $\la \in (0,1]$.
\end{abstract}

\renewcommand{\abstractname}{R\'esum\'e}
\begin{abstract}
On \'etudie la dynamique des produits crois\'es agissant sur le cylindre $\cyl$, de la forme
$$
\tht \mapsto (\ell \theta, \la \theta + \tau ( \theta)),
$$
o\`u $\ell \ge 2$ est un entier, $\la \in (0, 1)$ et $\tau: \T \to \R$ est une fonction continue.
On s'int\'eresse aux propri\'et\'es {\it topologiques} de l'attracteur global $\Omegalt$ de cet endomorphisme.
\'Etant donn\'e $\ell$ et une fonction lipschitzienne $\tau$, on d\'emontre que l'attracteur $\Omegalt$ est hom\'eomorphe \`a un anneau topologique pour tout $\lambda$ assez proche de~$1$.
D'autre part, on d\'emontre qu'il existe au plus un nombre fini de $\lambda \in (0, 1)$ tels que l'attracteur $\Omegalt$ soit une courbe de Jordan.

Ces r\'esultats s'appuient sur une analyse d\'etaill\'ee des \'equations {``cohomologiques''} it\'er\'ees~: $\tau = \cL_{\lambda_1} \mu_1$, $\mu_1 = \cL_{\lambda_2} \mu_2, \ldots$, o\`u $\cL_\lambda \mu = \mu \circ \ml - \lambda \mu$ et $\ml$ est l'application de multiplication par $\ell$ sur le cercle $\T$.
On d\'emontre le r\'esultat de finitude suivant~: toute fonction lipschitizenne $\tau$ s'\'ecrit de fa\c{c}on canonique sous la forme
$$
\tau = \cL_{\lambda_1} \circ \cdots \circ \cL_{\lambda_m} \mu,
$$
o\`u $m \ge 0$, $\lambda_1, \ldots, \lambda_m \in (0, 1]$ et la fonction lispchitzienne $\mu$ satisfait $\mu \neq \cL_\lambda \rho$ pour toute fonction continue $\rho$ et tout $\lambda \in (0, 1]$. 
\end{abstract}

\setcounter{tocdepth}{1}


\section{Introduction.}\label{affine-s}
In this paper we study the dynamics of skew product endomorphisms of the cylinder $\cyl$ of the form
$$
\begin{array}{rccc}
\Alt :& \cyl & \rightarrow & \cyl \\
& (\theta, t) & \mapsto &  ( \ell \theta , \la t + \tau ( \theta)),
\end{array}
$$
where $\ell \ge 2$ is an integer, $\la \in (0, 1)$ and  $\tau : \T \rightarrow \R$ is a continuous function.

The non-wandering set $\Omegalt$ of $\Alt$ is a global attractor of the dynamics of $\Alt$: the forward orbit of every point in $\cyl$ converges to $\Omegalt$ and $\Alt$ is transitive on $\Omegalt$.
In fact  $\Alt$ is topologically semi-conjugate to a solenoidal map on $\Omegalt$ (Section~\ref{preliminaries-s})

These maps where intially studied in~\cite{tsujii-01}, from a {\it measure theoretical} point of view.
In that paper M.~Tsujii showed that $\Alt$ has a unique physical measure and that the support of this measure is the attractor $\Omegalt$.
The main result of~\cite{tsujii-01} is that, when $\lambda > \ell^{-1}$, for generic functions $\tau$ of class ${\mathcal C}^2$ the unique physical measure of $\Alt$ is absolutely continuous with respect to Lebesgue measure.

The purpose of this paper is to study {\it topological} properties of the attractor sets $\Omegalt$.
Our main result is the following.

\begin{theorem}
\label{affine-th}
Suppose that $\tau: \T \rightarrow \R$ is Lipschitz. Then the following hold:

{\rm (1)} The set $\mathcal{J}_\tau = \{ \lambda \in (0,1) / \Omegalt \mbox{ is a Jordan curve} \}$ is finite.

{\rm (2)} There exists $\lambdanot \in (0,1)$ such that $\Omegalt$ is a homeomorphis to a closed  annulus for
all $\la \in [\lambdanot , 1)$.  
\end{theorem}

For a given $ \la \in (0,1)$, we characterize those functions $\tau$ for which $\Omegalt$ is a Jordan curve  in terms of the Fourier coefficients of $\tau$ (Theorem~2).
>From this characterization it follows easily that the set of those $\tau$ for which $\Omegalt$ is a Jordan curve has infinite codimension in the space of all Lipschitz functions.
\subsection{On the interior of the attractor.}
When $\lambda < \ell^{-1}$, it is easy to see that $\Omegalt$ has zero  Lebesgue measure and hence empty interior.
On the other hand, when $\lambda = \ell^{-1}$ we show the set $\Omegalt$ is not homeomorphic to an annulus (Proposition~\ref{ell-p}).
So $\lambdanot$ in part~(2) of the theorem must be strictly larger than~$\ell^{-1}$.

When $\lambda > \ell^{-1}$, Tsujii's result (mentioned above) implies that for most $\tau$ of class ${\mathcal C}^2$, the attractor  set $\Omegalt$ has positive Lebesgue measure.
Here we show examples of maps $\Alt$ for which the set $\Omegalt$ has non-empty interior but it is not homeomorphic to an annulus (Proposition~\ref{fat-hole-ex}).
In these examples $\lambda$ can be taken arbitrarily close to $\ell^{-1}$. 
Moreover we show that these examples are {\it robust} in the sense that any map $A : \cyl \to \cyl$ that is sufficiently (Lipschitz) close to $\Alt$ has the same properties.

In a forthcoming paper we show that, when $\lambda > \ell^{-1}$, for most $\tau$ of class ${\mathcal C}^2$ the set $\Omegalt$ has non-empty interior.
\subsection{On the iterated cohomological equation.}
\label{iter-ss}
Recall that a continuous function $\tau : \T \to \R$ is {\it cohomologous to}~0 if there is a continuous function $\mu$ such that
$$
\tau = \cL \mu = \mu - \mu \circ \ml,
$$
where $\ml : \T \to \T$ is the multiplication by $\ell$ map.
It is easy to see that in that case $\int_{\T} \tau = 0$ and the function $\mu$ is unique up to an additive constant.
For this reason we will assume that all the functions considered in the rest of the introduction have~0 integral.

Part~(2) of the theorem is first proven in the case when $\tau$ is not {\it cohomologous to}~0 (Proposition~\ref{not-coboundary-p}).
When $\tau = \cL \mu$ is cohomologous to~0, a direct computation shows that {\it the maps $\Alt$ and $A_{\lambda, \mu}$ are conjugated} by the homeomorphism $(\theta, t) \mapsto (\theta, (t + \mu)/(1 - \lambda))$.
So, if $\mu$ is not cohomologous to~0 we reduce to the first case.
By induction, if for some positive integer $n$ there is a continuous function $\mu : \T \longrightarrow \R$ that is not cohomologous to~0 and such that $\tau = \cL^n \mu$, then we reduce to the first case.

We complete the proof of part~(2) of the theorem  by showing that {\it a non-constant Lipschitz function cannot be infinitely cohomologous to~0}.
More precisely, we show that if $\tau$ is Lipschitz, then the integer $n$ above is bounded by a constant depending only on~$\tau$ (Lemma~\ref{finite-boundary-l}).

\

{\par\noindent {\bf Problem 1.}}
{\it Is there a non-constant continuous function that is infinitely cohomologous to~0?}
\subsection{Cohomological operators.}
\label{coho-ss}
For $\lainzoc$ it is interesting to consider the linear operators $\cL_\lambda$ defined by $\cL_\lambda \mu = \lambda \mu - \mu \circ \ml$, so that $\cL_1 = \cL$.
For $\la_0 \in (0,1)$ we show that $\Omega_{\lambda_0, \tau}$ is a Jordan curve if and only if there exists a (Lipschitz) continuous function $\mu$ such that $\tau = \cL_{\lambda_0} \mu$ (Proposition~\ref{jordan-characterization-p}).
In that case, for every $\lainzoo$ different from $\lambda_0$, the maps $A_{\lambda_0, \tau}$ and $\Alt$ are conjugate (Lemma~\ref{affine-conjugacy-l}).

We show that each Lipschitz function $\tau$ can be written in a canonical way as
$$
\tau = \cL_{\lambda_1} \circ \ldots \circ \cL_{\lambda_m} \mu,
$$
where the function $\mu$ satisfies $\int_{\T} \mu = 0$ and $\mu \neq \cL_\lambda \rho$ for every $\lainzoc$ and every continuous function $\rho$ (here there might be repetitions among $\lambda_1, \ldots, \lambda_m$).
This implies part~(1) of the theorem.
Note that such a function $\rho$ is such that for every $\lainzoo$ the set $\Omega_{\lambda, \mu}$ is not a Jordan curve and for every $\lainzoo$ different from the $\lambda_i$, the maps $\Alt$ and $A_{\lambda, \tau}$ are conjugate (Theorem~3).

\

{\par\noindent {\bf Problem 2.}}
{\it For $\ell =2$, let $\mu$ be a Lipschitz function such that $\mu \neq \cL_{\lambda} \rho$ for every $\lainzoo$ and every non-constant function $\rho$.
Is there $\la_0 \in (0,1)$ such that $\Omega_{\lambda, \mu}$ is homeomorphic to an annulus if and only if $\lambda \ge \lambda_0$?}

\subsection{Are there periodic points in the boundary?}
The {\it upper} (resp. {\it lower}) {\it boundary of the attractor} is by definition the graph of the function
$$
\rho^+(\theta) = \sup \{ t / (\theta, t) \in \Omegalt \}
$$
$$
(\mbox{resp. } \rho^- (\theta) = \inf \{ t / (\theta, t) \in \Omegalt \}.)
$$
These functions are continuous and characterized by the functional equations
$$
\rho^+(\theta) = \max \{ \lambda \rho(\theta') + \tau(\theta') / \theta' \in \ml^{-1}(\theta) \}
$$
$$
\rho^-(\theta) = \min \{ \lambda \rho(\theta') + \tau(\theta') / \theta' \in \ml^{-1}(\theta) \}.
$$
Moreover, when $\tau$ is Lipschitz the functions $\rho^+$ and $\rho^-$ are also Lipschitz.

\

{\par\noindent {\bf Conjecture.}}
{\it For each $\lambda \in (0, 1)$ there is an open and dense set of functions $\tau$ of class $\cC^1$, such that the following properties hold.
\begin{enumerate}
\item[1.]
The upper (resp. lower) boundary contains a finite number of periodic orbits of $\Alt$.
\item[2.]
The upper (resp. lower) boundary is formed by a finite number of pieces of the unstable manifolds of the periodic orbits that it contains.
\end{enumerate}
In particular the upper and lower boundaries are $\cC^1$ by parts.}

\

For a given $\lambda \in (0, 1)$ and a continuous function $\tau : \T \to \R$ consider the closed set
$$
D^+ = \{ \theta \in \T  /  \lambda \rho^+(\theta) + \tau(\theta) = \rho^+(\ell \theta) \}.
$$
It follows from the functional equation of $\rho^+$ that $\ml(D^+) = \T$, so the maximal invariant set
$$
K^+ = \{ \theta \in D^+  /  \ml^n(\theta) \in D^+ \mbox{ for } n \ge 1 \}
$$
is non-empty and compact.
Part~1 of the conjecture implies that $K^+$ contains a finite number of periodic orbits and part~2 implies that $K^+$ is finite.

The above conjecture is somewhat similar to the conjecture that, for generic expanding endomorphisms of the circle $f$, there is a unique measure $\mu$ minimizing the integral $\int_\T \ln f' d \mu$ and that this measure has finite support, see~\cite{bousch-00}, ~\cite{contreras-01} and references therein.  

\subsection{Notes and references.}
Similar skew product of the cylinder where studied by M.~Viana~\cite{viana-97}.
M. Tsujii~\cite{tsujii-03} extended the results of~\cite{tsujii-01}  and~\cite{viana-97} to general partially hyperbolic endomorphisms on surfaces.
Piecewise affine endomorphisms having an attractor of non-empty interior were studied in~\cite{dobrynski-99}.

\medskip
\noindent
{\bf Acknowledgements.}
{Juan Rivera-Letelier is grateful to J.-C. Yoccoz for useful conversations and valuable comments and Rodrigo Bam\'on thanks Marcelo Viana for helpful and constructive discussions.}

\subsection{Outline.}
Let us now describe the structure of the paper.

\medskip
Section~\ref{preliminaries-s} starts giving several equivalent characterizations of the 
set $\Omegalt$. Then, in Subsection~\ref{semi-ss}, we show that the semiconjugacy 
between the dynamics induced by the multiplication by $\ell$ map on
the solenoid $\cS$ and the dynamics of $\Alt$ on $\Omegalt$ can be written in a fairly explicit manner.
Also, we endow  $\cS$ with an adapted metric $\dist_\la$ which makes this semiconjugacy a Lipschitz function,
provided that $\tau$ is Lipschitz. Then it naturally follows 
 that the upper and lower boundaries of
the attractor $\Omegalt$ are Lipschitz graphs when $\tau$ is Lipschitz  (Subsection~\ref{uplow-ss}).

\medskip
Section~\ref{jordan-s} is devoted to study Jordan curve attractors. We characterize them
and show, among other results,  that $\Omegalt$ is a Jordan curve if and only if the
functional equation $\mu \circ \ml - \la \mu = \tau$ has a continuous solution (Proposition~\ref{jordan-characterization-p}). 
This allow us to show, in Subsection~\ref{codim-ss}, that the set of continuous $\tau$ (with absolutely convergent Fourier series)
such that $\Omegalt$ is a Jordan curve has infinite codimension (Theorem~\ref{codim-th}).

\medskip
The main result in Section~\ref{annulus-s}  is that, for $\la$ sufficiently close to $1$,
the attractor $\Omegalt$ is a closed topological annulus provided that $\tau$ is H\"older, not cohomologous to $0$ 
and with $0$ integral (Proposition~\ref{not-coboundary-p}). 
This section starts with general results about annular attractors.
In particular, we show that if the image of the upper boundary of $\Omegalt$ is 
above  in the cylinder $\cyl$ than the image of the lower boundary, then $\Omegalt$ is homeomorphic
to a closed annulus (Lemma~\ref{annulus-characterization-l}). 
Also, we establish that for $\la \leq 1/\ell$, the attractor $\Omegalt$ cannot
be an annulus (Proposition~\ref{ell-p}). Then, in Subsection~\ref{not-coboundary-ss}, 
under the above assumptions for  $\tau$  we find periodic orbits $\cO^\pm$ in the circle so that the corresponding
orbits in $\Omegalt \subset \cyl$ have $\R$ coordinates tending to $\pm \infty$ as $\la \rightarrow 1$.
>From this we deduce that the image of the upper boundary is higher up in $\cyl$ than the image of 
the lower boundary, when $\la$ is sufficiently close to $1$, and therefore that $\Omegalt$ is an annulus.

\medskip
Sections~\ref{jordan-s} and~\ref{annulus-s} lead us to study in more detail the operators
$\cL_\la \mu = \mu \circ \ml - \la \mu$ with $\la \in (0,1]$. For our purpose the natural domain
of the operators $\cL_\la$ is the space of Lipschitz functions.
In Section~\ref{linear-s} we start by proving some general facts about these linear operators
and relating them to conjugacy classes of maps of the form $\Alt$ (lemmas~\ref{lla-l} and~\ref{affine-conjugacy-l}).
As mentioned in Subsection~\ref{iter-ss}  this forces us to study iterated equations of the
form $\cL_{\la_1} \mu_1 = \tau$, $\cL_{\la_2} \mu_2 = \mu_1$, $\dots$.
An important feature of the operators $\cL_\la$ is that  they do  not 
increase the best Lipschitz constant  for
$\mu$ (Lemma~\ref{lip-norm-l}) which implies that  solutions $\mu_n$ of the  iterated equations
are uniformly Lipschitz. 
Then we show that the above equations have the effect of increasing the Fourier coefficients of $\mu_n$ as 
$n$ increases and  establish our Main Lemma  which states that, given a Lipschitz function $\tau$
there exists a finite collection $0 < \la_1, \dots, \la_m \leq 1$ (maybe with repetitions) and a (Lipschitz)
continuous function $\mu$ such that $\cL_{\la_1} \circ \cdots \circ  \cL_{\la_m} \mu = \tau$ and 
$\cL_\la \rho \neq \mu$ for all continuous $\rho$ (see lemmas~\ref{finite-boundary-l} and~\ref{finite-jordan-l}).
Theorem~\ref{affine-th} and its stronger version Theorem~\ref{strong-affine-th} follow immediately
from our main lemma.

\medskip
In Section~\ref{uplow-s} we start by appropriately defining the attractor set  and the upper and lower boundaries
for maps which are close to $\Alt$
and state that the upper and lower boundaries of $\Omegalt$ vary continuously
in the $\cC^0$ topology under Lipschitz perturbations of $\Alt$ (Proposition~\ref{uplow-p}).
To prove this we pass to the universal cover $\R^2$  of $\cyl$ and, in Subsection~\ref{r2-ss}, examine the action of
Lipschitz maps from $\R^2$ into $\R^2$ on the graphs of Lipschitz functions. 
In Subsection~\ref{fixed-ss},
motivated by the fact that the upper and lower boundaries of $\Omegalt$ are the graphs of
functions $\rho^\pm : \T \rightarrow \R$ which satisfy certain functional equations we
show that, under certain conditions, the upper and lower boundaries 
of the attractor of maps $F$ of the cylinder are fixed points
of operators $\cT^\pm_F$ which act on Lipschitz functions. 
The definition and properties of $\cT^\pm_F $ rely on the work of Subsection~\ref{r2-ss}.
At the end of  Section~\ref{uplow-s} we prove  the above mentioned  continuity of the upper and lower 
boundaries of the attractor.

\medskip
The last section, Section~\ref{example-s}, contains two examples. 
The first example consists of an application of our results to the study of
a family $f_{\la, c}$ of endomorphisms of $\Cstar = \C \setminus \{ 0 \}$ where $\la \in (0,1)$ and
$c \in \Cstar$. Here $f_{\la, c} = f_{\la, 0} + c$ where $f_{\la,0}$ acts as 
angle doubling on the arguments of $z \in \Cstar$ and as an  affine contraction of factor $\la$ in the radial direction.
Thus, $f_{\la, c}$ is closely related to the extensively studied quadratic family $Q_c(z)= z^2 + c$ where 
the $|z| \mapsto |z|^2$ action of $Q_0$ in the radial direction has been  substituted by an affine contraction. We show 
that for $\la$ sufficiently close to $1$ and for $|c|$ sufficiently small, 
the attractor of $f_{\la, c}$ is an annulus. The second example shows that given $\la > 1/\ell$
there exist skew product endomorphisms $\Alt$ such that the attractors $\Omegalt$ have non-empty interior and are 
not an annulus.  This example is robust under Lipschitz perturbations. 


\section{Preliminaries}
\label{preliminaries-s}
Throughout this section, unless otherwise stated, $\tau : \T \rightarrow \R$ is a continuous
function and $\la \in (0,1)$. 
We start by showing that $\Omegalt$ is a global attractor for the dynamics of 
$\Alt$ and giving several equivalent characterizations of this set.
As usual $\| \tau \|_\infty = \sup \{ \tau (\theta) \,\, / \, \theta \in \T \}$.

\begin{lemma}
  \label{global-attractor-l}
Let $\Omegalt$ be the non-wandering set of $\Alt : \cyl \rightarrow \cyl$.

{\rm (1)} 
If   $U_0 = \T \times (-T_0,T_0)$  for some $T_0$ such that $(1-\la) T_0 > \|\tau\|_\infty$, 
then $$\Alt(\overline{U_0}) \subset U_0 \mbox{  and } \Omegalt = \cap_{n \geq 0}  \Alt (U_0).$$

{\rm (2)} $\Omegalt$ is the set of all $\tht \in \cyl$ with a bounded infinite backward orbit (i.e.,
there exists $C >0$ and $\{ (\theta_n,t_n) \}_{n \geq 1}$ such that $\Alt^n (\theta_n,t_n)= (\theta,t)$
and $|t_n| \leq C$ for all $n \geq 1$).

{\rm (3)} $\Omegalt$ is the closure of the set formed by the periodic points of $\Alt$.
\end{lemma}

\noindent
{\bf Proof.} 
Denote by $\operatorname{Per}$ the set of periodic points of $\Alt$ and by $B$ the set of points in $\cyl$ which
have a bounded infinite backward orbit. We will show that:
$$\overline{\operatorname{Per}} \subset \Omegalt \subset \cap_{n \geq 0} A^n_{\la, \tau} ({U_0})  \subset B \subset \overline{\operatorname{Per}}.$$
The inclusions $\overline{\operatorname{Per}} \subset \Omegalt $ and $ \cap_{n \geq 0} A^n_{\la, \tau} ({U_0})  \subset B$ are clear.

Let  $T_0$ be  such that $(1-\la) T_0 > \| \tau \|_\infty$. Since
$$|\la t + \tau (\theta)| \leq \la T_0 + \| \tau \|_\infty < \la T_0 + (1-\la) T_0 = T_0.$$
for all $|t| \leq T_0$, it follows that  $\Alt(\overline{U_0}) \subset U_0$.

Note that if $(\theta_n, t_n) = A^n_{\la,\tau} (\theta, t)$, then $\theta_n = \ell^n \theta$ and
$$  t_n  =  \la^n t + \la^{n-1} \tau (\theta) + \la^{n-2} \tau (\ell \theta) + \cdots +  \tau (\ell^{n-1} \theta). $$
Therefore, 
\begin{equation}
  |t_n|  \leq   \frac{1-\la^n}{1-\la} \| \tau \|_\infty + \la^n |t|. 
\label{global-attractor-eqn}
\end{equation}
Hence, for all $\tht \in \cyl$ there exists $n$ such that $A^n_{\la, \tau} \tht \in U_0$.

Now we show that $\Omegalt \subset \cap A^n_{\la, \tau} (U_0)$. In
fact, suppose that $(\theta,t) \notin A^m_{\la, \tau} (U_0)$ for some
$m$ and consider $n$ such that $\tht$ belongs to the open set $V=
\Alt^{-n} (U_0)$. It follows that $\tht$ does not belong to the
non-wandering set $\Omegalt$ since $\tht \in V$ but $\tht \notin
\Alt^{n+M} (V)$ for all $M \geq m$.

To finish  the proof of the lemma we consider a neighborhood $U$ of a point $\tht$ with bounded infinite backward orbit
and proceed to show that $U$ contains a periodic point.
Let $\{ (\theta_n,t_n) \}_{n \geq 1}$ and $C>0$ be  such that $\Alt^n (\theta_n,t_n)= (\theta,t)$
and $|t_n| \leq C$ for all $n \geq 1$. 
There exist  $n \geq 1$ and an open interval $I \subset \T$ around $\theta_n$ such
that the rectangle $R= I \times [-T_0, T_0]$ maps into $U$ under $A^n_{\la,\tau}$. Since there exists
$\theta^\prime \in I$ periodic under multiplication by $\ell$, say of period $m$, we have that $A^m_{\la, \tau}$ restricted to 
$\{\theta^\prime\}  \times [-T_0,T_0]$ is a contraction. It follows that  $\Alt$ has a periodic point in 
$\{\theta^\prime\}  \times [-T_0,T_0] \subset R$ 
and therefore in $U$.
\hfill $\Box$

\subsection{The solenoid and the semiconjugacy.}
\label{semi-ss}
Throughout the paper multiplication by $\ell$ in the circle will be denoted by $\m_\ell: \T \rightarrow \T$. 
For each $\lainzoo$ we endow the solenoid:
$$\cS := \{ \bar{\theta} = (\theta_k) \in (\T)^{\N \cup \{ 0 \}} \, / \m_\ell( \theta_{k+1} )= \theta_k \mbox{ for all } k \geq 0 \}$$
with the  adapted metric
$$\ladist ( (\theta_k), (\thetap_k)) = \sum_{k \geq 0} \la^k \Tdist (\theta_k, \thetap_k).$$
where $\Tdist$ denotes the projection of the standard metric of $\R$ onto $\T$.
The dynamics of multiplication by $\ell$ induces:
$$\begin{array}{rclc}
\cM_\ell :& \cS & \rightarrow & \cS \\
& (\theta_k)_{k \geq 0} & \mapsto &  (\ell \thetanot, \thetanot, \theta_1, \dots ).
\end{array}$$

We will show that the dynamics of  $\cM_\ell: \cS \rightarrow \cS$  semi-conjugates to that of
$\Alt : \Omegalt \rightarrow \Omegalt$. Thus, the attractor $\Omegalt$ is in this sense a solenoidal attractor. 
To write an explicit formula for the semiconjugacy from $\cS$ onto $\Alt$ we need the following definition.

\begin{definition}
  Given  a continuous function  $\tau: \T \rightarrow \R$ and  $\lainzoo$ we define
$t_\la : \cS \rightarrow \R$ by
$$t_\la (\btheta) = \tau(\theta_1) + \la \tau (\theta_2) + \la^2 \tau (\theta_3) + \cdots. $$
\end{definition}

Note that $t_\la$ is continuous. Under the assumption that $\tau$ is a Lipschitz function 
we will show that $t_\la$ is also Lipschitz. In order to make the statements precise we introduce
some notation.

\begin{notation}
  {\em Consider two metric spaces $(X,\rho_X)$, $(Y,\rho_Y)$.
Given a Lipschitz map $f: X \rightarrow Y$ the best Lipschitz constant for $f$ 
$$\sup_{a \neq b} \frac{\rho_Y(f(a),f(b))}{\rho_{X}(a,b)}$$
is denoted by $\|f\|_L$ and if $C \geq \| f \|_L$, we say that $f$ is a $C$--Lipschitz map.}
\end{notation}

\begin{lemma}
\label{tla-lip-l}
If  $\lainzoo$ and $\tau : \T \rightarrow \R$ is  a $\| \tau \|_L$--Lipschitz map, 
then $t_\la$ is a  $(\la^{-1} \| \tau \|_L) $--Lipschitz from $(\cS, \ladist)$  to $\R$.
\end{lemma}

\proof
For any $\btheta = (\theta_k)$ and $\btheta^\prime = (\thetap_k)$ in $\cS$ we have
\begin{eqnarray}
|t_\la (\btheta) - t_\la (\btheta^\prime)| & \leq & \sum_{k \geq 1}  \la^{k-1} | \tau (\theta_k) - \tau(\thetap_k)| \nonumber \\
                                           & \leq & \| \tau \|_L \cdot \sum_{k \geq 1} \la^{k-1} \Tdist (\theta_k,\thetap_k) \nonumber \\
                                           & = & \la^{-1} \| \tau \|_L (\ladist(\btheta, \btheta^\prime) - \Tdist (\theta_0, \thetap_0) \label{tla-lip-e} \\
& \leq & \la^{-1} \| \tau \|_L \ladist(\btheta, \btheta^\prime) \nonumber
\end{eqnarray}
That is, $t_\la$ is $\la^{-1} \| \tau \|_L$--Lipschitz.
\hfill $\Box$

\begin{proposition}
  \label{affine-semiconjugacy-p}
Given $\lainzoo$ and a continuous function  $\tau : \T \rightarrow \R$,  let 
$$\begin{array}{rclc}
h_\la :& \cS & \rightarrow & \cyl \\
& \btheta= (\theta_k) & \mapsto &  (\thetanot, t_\la (\btheta)).
\end{array}$$
Then $h_\la$ is a continuous semiconjugacy  from $\cS$ onto $\Omegalt$.
That is, $h_\la: \cS \rightarrow \Omegalt$ is onto and $\Alt \circ h_\la = h_\la \circ \cM_\ell$.
Moreover, $h_\la : (\cS, \dist_\la) \rightarrow \Omegalt$ is Lipschitz whenever 
$\tau : \T \rightarrow \R$ is Lipschitz. 
\end{proposition}

\proof 
That $\Alt \circ h_\la = h_\la \circ \cM_\ell$ is a straightforward
computation.  We must show that $h_\la (\cS) =
\Omegalt$.  Since $h_\la (\cS) $ is bounded and forward invariant (i.e., $\Alt(h_\la (\cS)) = h_\la (\cS)$),
by Lemma~\ref{global-attractor-l}, we have that $h(\cS) \subset
\Omegalt$.  Now if $(\thetanot, \tnot) \in \Omegalt$, then there exists
a bounded backward orbit $\{ (\theta_n, t_n)\}_{n \geq 0}$. Therefore,
$\btheta = (\theta_n) \in \cS$ and $\tnot = t_\la (\btheta)$.  Hence,
$(\thetanot, \tnot) = h_\la (\btheta) \in h_\la (\cS)$.  
By the previous lemma, $h_\la$ is Lipschitz, if $\tau :\T \rightarrow \R$
is a Lipschitz function. 
\hfill $\Box$

\subsection{Upper and lower boundaries of the attractor.}
\label{uplow-ss}
The attractor $\Omegalt$ lies in between the graph of two functions which we call the upper and lower boundaries
of $\Omegalt$. More precisely:

\begin{definition}
Let 
\begin{eqnarray*}
  \rho^+ (\theta) &=& \sup\{ t \, / \, \tht \in \Omegalt \}, \\
  \rho^- (\theta) &=& \inf \{ t \, / \, \tht \in \Omegalt \}. \\ 
\end{eqnarray*}
  We say that $\partial_{\pm} \Omegalt =\{ (\theta, \rho_{\pm}(\theta) \, / \,  \theta \in \T \}$
are the {\bf upper and lower boundaries of $\Omegalt$}, respectively.
\end{definition}

Since $\Alt (\Omegalt) = \Omegalt$ and $\Alt$ is locally orientation preserving, it follows that
\begin{eqnarray}
\label{functional-e}
 \rho^+ (\theta) & =& \max \{ \la \rho^+ (\thetap) + \tau(\thetap) \,\,/\, \thetap \in \m^{-1}_\ell (\theta) \}, \\
 \rho^- (\theta) & =& \min \{ \la \rho^- (\thetap) + \tau(\thetap) \,\,/\, \thetap \in \m^{-1}_\ell (\theta) \}. \nonumber
\end{eqnarray}
It is interesting to observe that the above functional equation is very similar to the one contained in  Lemma A of~\cite{bousch-00}. 

\smallskip
 For $\tau$ Lipschitz, the upper and lower boundaries are Lipschitz graphs. In general, for $\tau$ of class $\cC^\infty$ or even
real analytic, the upper and lower boundaries are not $\cC^1$. 

\begin{lemma}
  \label{affine-upper-lower-l}
If $\tau : \T \rightarrow \R$ is Lipschitz,
then $\rho^\pm : \T \rightarrow \R$ are $
\frac{\| \tau \|_L}{(\ell-\la)}$--Lipschitz maps.
\end{lemma}

\proof
Consider  $\theta_0, \thetap_0 \in \T$.
Let $\btheta = (\theta_k) \in \cS$ be such that $t_\la (\btheta)= \rho^+(\thetanot)$.
There exists $\bthetap = (\thetap_k) \in \cS$
such that  $\Tdist (\thetap_k, \theta_k) = \ell^{-k} \Tdist (\thetap_0, \theta_0)$ and therefore, 
$$\ladist (\btheta, \bthetap) = \frac{\ell}{\ell-\la} \Tdist(\thetanot, \thetap_0).$$
Then, by (\ref{tla-lip-e}):
\begin{eqnarray*}
\rho^+ (\thetap_0) \geq t_\la (\btheta^\prime) & = & t_\la (\btheta^\prime) -  t_\la (\btheta) +  t_\la (\btheta) \\
& \geq & -C (\ladist(\btheta,\btheta^\prime) - \Tdist(\thetanot, \thetap_0)) + t_\la (\btheta)\\
& = & -C (\frac{\ell}{\ell-\la} -1)\ladist(\btheta,\btheta^\prime) + t_\la (\btheta)\\ 
& =  &\rho^+ (\thetanot)  -C \frac{\la}{\ell-\la} \ladist(\btheta,\btheta^\prime)
\end{eqnarray*}
where $C= \la^{-1} \| \tau \|_L $ is a Lipschitz constant for $t_\la$ 
(see Lemma~\ref{tla-lip-l}).
It follows that $\rho^+ : \T \rightarrow \R$ is Lipschitz
with the appropriate constant. For  $\rho^- : \T \rightarrow \R$
a similar argument can be applied.
\hfill $\Box$

\section{Jordan curve attractors}
\label{jordan-s}

Observe that the equator $\{0 \} \times \T$ of the cylinder
 is the attractor of the map  $A_{\la, 0}$.

\subsection{Characterization}
The next proposition characterizes Jordan curve attractors.
\begin{proposition}
\label{jordan-characterization-p}
Let $\tau : \T \rightarrow \R$ be a function of class $\cC$
where $\cC$ is either the Lipschitz class, or $\cC^r$ class for
some $r \in [0, \infty] \cup \{\omega\}$.

Then the following are equivalent:

{\rm (1)} $\Alt$ is topologically conjugate to $A_{\la,0}$.

{\rm (1')}   $\Alt$ is $\cC$-conjugate to $A_{\la,0}$.

{\rm (2)} $\Omegalt \subset \cyl$ is the graph of a continuous function $\mu: \T \rightarrow \R$.

{\rm (2')} $\Omegalt \subset \cyl$ is the graph of a function $\mu: \T \rightarrow \R$ of class $\cC$.

{\rm (3)} The functional equation $\mu \circ \m_\ell - \la \mu = \tau$ has a continuous solution $\mu: \T \rightarrow \R$.

{\rm (3')} The functional equation $\mu \circ \m_\ell - \la \mu = \tau$ has a  solution $\mu: \T \rightarrow \R$ of class $\cC$.

{\rm (4)} $\rho^+(\theta) = \rho^- (\theta)$ for all  $\theta \in \T$.

{\rm (5)} $\rho^+ (\theta) = \rho^- (\theta)$ for some $\theta \in \T$.

{\rm (6)} $\Omegalt$ is a Jordan curve.

\end{proposition}

\proof
(6)$\implies$(5). If $\Omegalt$ is a Jordan curve, then $\rho^+ (\theta) = \rho^- (\theta)$ for some $\theta \in \T$.
Otherwise, $\rho^+(\theta) > \rho^-(\theta)$ for all $\theta$ and  the two Jordan
curves $\{ (\theta, \rho^+(\theta)) \}$ and $\{ (\theta, \rho^-(\theta)) \}$ would be disjoint and contained
in the  Jordan curve $\Omegalt$ which is impossible.

(5)$\implies$(4). If $\rho^+ (\theta) = \rho^- (\theta)$ for some $\theta \in \T$, then  the set
$S= \{ \theta \in \T \, /\, \rho^+ (\theta) = \rho^- (\theta) \}$ is not empty, closed
and $\m^{-1}_\ell (S) = S$, hence  $S = \T$.

(4)$\implies$(2). If $\rho^+$ and  $\rho^-$ agree on $\T$, then $\Omegalt$ is the graph of $\mu = \rho^+ = \rho^-$.

(2)$\implies$(3) (resp. (2')$\implies$(3')). 
If $\Omegalt \subset \cyl$ is the graph of a $\cC^0$ (resp. $\cC$) function $\mu: \T \rightarrow \R$,
then $(\ell \theta , \la \mu (\theta) + \tau (\theta)) = \Alt(\theta, \mu(\theta))$ belongs to the graph of $\mu$.
Therefore, $\mu (\ell \theta) = \la \mu (\theta) + \tau (\theta)$ for all $\theta$.

(3)$\implies$(3'). Let us denote by  $\pi : \R \rightarrow \T$ quotient map. If
$\mu \circ \m_\ell - \la \mu = \tau$ has a $\cC^0$ solution $\mu: \T \rightarrow \R$,
then $\mu \circ \pi: \R \rightarrow \R$ is a solution of $\tmu (\ell s) - \la \tmu(s) = \ttau(s)$
where $\ttau =  \tau \circ \pi $. It is not difficult to 
check that this latter equation has a unique continuous solution   given by
$$\tmu (s) = \ttau(s/\ell) + \la \ttau(s/\ell^2) + \la^2 \ttau(s/\ell^3) + \cdots,$$
which is of the same class as $\ttau$. Therefore, $\mu \circ \pi = \tmu$ and $\mu$ are of class $\cC$.

(3')$\implies$(1') and (2'). If $\mu: \T \rightarrow \R$ is a $\cC$ function such that  $\mu \circ \m_\ell - \la \mu = \tau$,
then $h \circ A_{\la, 0} = \Alt \circ h$ where $h(\theta,t) = (\theta, t + \mu (\theta))$. Hence $\Alt$ is 
$\cC$-conjugate  to $A_{\la,0}$ and $\Omegalt=h(\T \times \{0\} = \Omega_{\la,0})$ is the graph
of $\mu: \T \rightarrow \R$.

Since (1') trivially implies (1) and (1) implies (6), the proof of the proposition is complete.
\hfill $\Box$

\subsection{Infinite codimension.}
\label{codim-ss}
Throughout this subsection we fix $\la \in (0,1)$. In 
Proposition~\ref{jordan-characterization-p} we showed  that the global  attractor $\Omegalt$ of $\Alt$ is a
Jordan curve if and only if the functional equation 
\begin{equation}
\tau  =\mu \circ \m_{\ell}
-\lambda \mu
\label{eqa1}
\end{equation}
has a continuous solution $\mu :\R/\Z\rightarrow \R$.
The aim of this subsection is to characterize the  continuous functions $\tau :
\R/\Z\rightarrow \R$ for which the functional equation~(\ref{eqa1}) has a solution.
With this purpose in mind, for each $k \in \N$, we introduce the linear functional

$$\begin{array}{rccl}
\cD_k :& L^1(\T,\R) & \rightarrow & \R \\
& \psi & \mapsto &  \int_\T \psi (\theta) {\nu_k(\theta)} d\theta
\end{array}$$
where 
$$ \nu_{k}( \theta ) =\sum_{n\geq 0}\lambda ^{n}( e^{2\pi i \ell^{n}k \theta}+e^{-2\pi i\ell^{n}k\theta}).$$

\begin{lemma}
  \label{measurable-l}
  Consider $\tau \in L^1(\T,\R)$ and suppose that there exists $\mu \in L^1(\T,\R)$
such that:
$$\tau  = \mu \circ \m_{\ell} -\lambda \mu.$$
Then, $\cD_k (\tau) =0$ for all $k \in \N$ such that $\ell \nmid \, k$.
\end{lemma}

Recall that for each $\phi \in L^{1}( \R/\Z, \R) $ the $k$-th Fourier coefficient of ${\phi}$
is defined by 
\begin{equation*}
\widehat{\phi}( k) =\int_{\T} {\phi}( \theta ) e^{-2\pi ik \theta}d\theta.
\end{equation*}
For general background on Fourier series see~\cite{katznelson-76}.

\proof
>From  equation $( \ref{eqa1})$ it follows that
\begin{equation}
\lambda ^{n}\tau ( \theta ) +\lambda^{n-1}\tau (\ell  \theta
) +\cdots  +\tau (
l^{n}\theta ) =\mu  (  \ell^{n+1}\theta ) -\lambda ^{n+1}\mu
 ( \theta )
\nonumber
\end{equation}
for all $n \geq 1$. For all $k$ such that
 $\ell \nmid k$, computing the $\ell^{n}k$-th Fourier coefficient of  the functions  involved in
the previous equation we obtain that:
\begin{equation}
\lambda ^{n}\hat{\tau} ( \ell^{n}k) +\lambda ^{n-1}\hat{\tau}
^{n-1} ( \ell^{n-1}k) +...+\hat{\tau} ( k) =-\lambda ^{n+1}
\hat{\mu} ( \ell^{n}k).
\label{eqa3}
\end{equation}
In view of the fact that    the Fourier coefficient of $\mu \in
L^{1} ( \R/\Z,\R)$ are bounded, as $n\rightarrow
\infty$ we have that $\lambda
^{n+1}\hat{\mu} ( \ell^{n}k) \rightarrow 0$ and, therefore, 
\begin{equation}
\sum_{n\geq 0}\lambda ^{n}\hat{\tau} ( \ell^{n}k) =0.  \label{eqa4}
\end{equation}

Now recall that the Fourier coefficients of a real valued function are even. 
In particular, $\hat{\tau} (k) =\hat{\tau} (-k)$ and equation $(\ref{eqa4})$ is equivalent to: 
\begin{equation}
\mathcal{D}_{k} ( \tau ) = \sum_{n\geq 0}\lambda ^{n} ( \hat{\tau} ( \ell^{n}k) +\hat{\tau}
 ( -\ell^{n}k) ) =0  
\nonumber
\end{equation}
Thus  we have completed the proof of the lemma.
\hfill $\Box$

\begin{theorem}
\label{codim-th}
Let $\tau :\R/\Z\rightarrow \R$ be a continuous
function with absolutely
convergent Fourier series. If  
\begin{equation}
\mathcal{D}_{k} ( \tau ) =0
\end{equation}
for all $k \in \N$ such that $\ell \nmid k$,
then there exists a continuous function $\mu :\R/\Z\rightarrow \R$ 
satisfying the functional equation $$\tau  = \mu  \circ \m_{\ell}  -\lambda \mu.$$
\end{theorem}

\proof
We will obtain $\mu$ as a Fourier series with coefficients $b(\cdot)$ where:

\begin{eqnarray}
b (0) & = & \dfrac{\hat{\tau} ( 0) }{1-\lambda }, \nonumber \\
b ( \ell^{n}k) & = & -\frac{1}{\lambda ^{n+1}}\sum_{j=0}^{n}\lambda ^{j} \hat{\tau} ( \ell^{j}k)  \label{eqa6}\\
 & = & \frac{1}{\lambda ^{n+1}}\sum_{j=n+1}^{\infty}\lambda ^{j}\hat{\tau} ( \ell^{j}k) \nonumber
\end{eqnarray}
for $n\geq 0$ and $k\neq 0$ such that $\ell$ does not divide $k$.
(Compare with equation (\ref{eqa3})).

We claim  that the Fourier series 
$$\sum_{q\in \Z} b ( q)
e^{2\pi iq \theta}=b ( 0) +\sum_{\ell\nmid k}\sum_{n\geq 0}b (
\ell^{n}k) e^{2\pi i \ell^{n} k\theta}
$$
is absolutely convergent (i.e.,  $\sum_{q\in \Z} |b (q) | <\infty$).
In fact, from  equation (\ref{eqa6}):
\begin{eqnarray*}
	 \sum_{q\in \Z} | b ( q) |& = &  | b ( 0) 
|+\sum_{\ell\nmid k}\sum_{n\geq 0} | b (\ell^{n}k) |
		 \\
	& \le & | b ( 0) | +\sum_{\ell\nmid k}\sum_{n\geq 0}\frac{1}{
\lambda ^{n+1}}\sum_{j\geq n+1}\lambda ^{j} | \hat{\tau} (
\ell^{j}k) |.
\end{eqnarray*}
Since 
\begin{eqnarray*}
\sum _{n\geq 0}\frac{1}{\lambda ^{n+1}}\sum _{j\geq n+1}\lambda
^j | \hat{\tau} ( \ell^j k) | & =& \sum _{n\geq 0} 
\sum _{j=0}^n \lambda ^j   | \hat{\tau} ( \ell^{n+1}k)
| \\ &=& \sum _{n\geq 0}\frac{1-\lambda ^{n+1}}{1-\lambda } |  \cdot \hat{\tau}
 ( \ell^{n+1}k) |\\ & \leq & \frac{1}{1-\lambda }
\sum _{n\geq 0} | \hat{\tau} ( \ell^{n+1}k) |,
\end{eqnarray*}
it follows that
\begin{eqnarray*}
\sum_{q\in \Z} | b ( q) | & \leq &  | \hat{\tau}
 ( 0) | \frac{1}{1-\lambda }+\frac{1}{1-\lambda }
\sum_{\ell\nmid k}\sum_{n\geq 0} | \hat{\tau} ( \ell^{n+1}k) | \\
& = & \frac{1}{1-\lambda }\sum_{q\in \Z} | \hat{\tau} ( q)
|  <\infty .
\end{eqnarray*}

Thus the series  $\sum_{q\in \Z}b ( q) e^{2\pi iq \theta}$
uniformly converges to a continuous function $\mu: \T \rightarrow \R$.
which  is a solution 
of $\tau =\mu \circ \m_{\ell}-\lambda \mu$ since an easy computation shows that
 the Fourier coefficients of $\tau $ and  $\mu \circ \m_{\ell}-\lambda \mu$ agree
(cf. Lemma~\ref{measurable-l}).
\hfill $\Box$

\begin{definition}
Let $\tau :\R/\Z\rightarrow \R$ be a continuous
function with  absolutely convergent Fourier series. 
Define the {\bf canonical representative of $\tau $}
by the continuous function $\tau _{c}:\R/\Z\rightarrow 
\R$ given by the Fourier
series ${\tau}_{c} ( \theta) =\sum_{\ell\nmid k}\mathcal{D}
_{k} ( \tau ) e^{2\pi ik\theta}$.
\end{definition}

An immediate consequence of the theorem above is the following corollary:
\begin{corollary}
Let  $\tau :\R/\Z \rightarrow \R$ be a continuous
function with  absolutely
convergent Fourier series. If $\tau _{c}$ is the canonical representative of 
$\tau $, then $\tau$ is cohomologous to  $\tau _{c}$. That is,  there exists  a
continuous function $\mu :\R/\Z\rightarrow \R$  such
that $\tau =\tau _{c}+\mu \circ m_{\ell}-\lambda \mu$.
\end{corollary}

\proof
Note that for every $k$ which is not divided by $\ell$ it holds that  $\mathcal{D}_{k} ( \tau -\tau _{c}) =0$. Therefore we may 
apply  the previous theorem to $\tau - \tau_c$ and obtain the corollary.
\hfill $\Box$

\section{Annular Attractors.}
\label{annulus-s}
Our next result contains the sufficient condition for $\Omegalt$ to be a closed topological annulus which will be used
in the proof of Theorem~\ref{affine-th}.  More precisely, below we prove that if the image of the upper boundary of $\Omegalt$ lies 
higher in $\cyl$ than the image of the lower boundary, then $\Omegalt$ is an annulus:
\begin{lemma}
  \label{annulus-characterization-l}
If $\Omegalt$ is not a Jordan curve and 
\begin{equation}
\min \{ \la \rho^+ (\thetap) + \tau (\thetap) \, / \, \thetap \in \m^{-1}_\ell (\theta) \} 
 \geq 
\max \{ \la \rho^- (\thetap) - \tau (\thetap) \, / \, \thetap \in \m^{-1}_\ell (\theta) \} 
\label{annular-e}
\end{equation}
for all $\theta \in \T$, then $\Omegalt$ is a closed topological annulus.
Moreover, if $\ell =2$ and $\Omegalt$ is a closed topological  annulus, 
then inequality (\ref{annular-e}) holds.
\end{lemma}

\begin{remark}
\label{annulus-r}
{\em For all $\ell \geq 2$ and provided that $\Omegalt$ is not a Jordan curve, a  necessary and sufficient condition for
$\Omegalt$ to be homeomorphic to a closed annulus is that for all $\theta \in \T$:
$$I (\theta)=\cup_{\thetap \in \m^{-1}_{\ell} (\theta)} \Alt( I(\thetap) )$$
where $I(\theta) = \{ \theta \} \times [\rho^- (\theta), \rho^+(\theta)]$.}
\end{remark}

\proof
Let $\cA = \{ \tht \, / \, \rho^-(\theta) \leq t \leq \rho^+ (\theta) \}$ and
note that $\Alt (\cA) \subset \cA$. 
We claim that from (\ref{annular-e}) it follows  $\Alt(\cA) = \cA$
and therefore $\Omegalt = \cA$. Otherwise there would exist
$t_0 \in (\rho^-(\theta_0), \rho^+(\theta_0))$ such that
$(\theta_0, t_0) \notin \Alt(\cA)$. Consider
$\thetap_{\pm} \in \m^{-1}_\ell (\theta_0)$ such that
$\rho^{\pm} (\theta_0) = \la \rho^{\pm}(\thetap_{\pm}) + \tau (\thetap_\pm)$.
Since $(\theta_0,t_0) \notin \Alt (\cA)$,
the image of $\cA \cap (\{ \thetap_+ \} \times \R)$ should be above 
$t_0$ and the image of $\cA \cap (\{ \thetap_- \} \times \R)$ should be 
below $t_0$. Hence, we would have:
$$ \la \rho^{+}(\thetap_{-}) + \tau (\thetap_-) < t_0 <  \la \rho^{-}(\thetap_{+}) + \tau (\thetap_+)$$
which contradicts (\ref{annular-e}). Therefore $\cA = \Omegalt$.
By the Proposition~\ref{jordan-characterization-p}, if $\Omegalt$ is not a Jordan curve, then
$\rho^+ > \rho^-$ on $\T$. Hence $\Omegalt = \cA$ is a closed topological  annulus.

For $\ell =2$, if (\ref{annular-e}) does not hold for some $\theta_0$, then
the two intervals, $\Alt(\cA \cap (\{\theta_0 /2 \} \times \R))$ and
$\Alt(\cA \cap (\{\theta_0 /2 +1/2 \} \times \R))$ are disjoint and
their union is $\Alt(\cA) \cap (\{ \theta_0 \} \times \R)$. Therefore,
$\Omegalt \subset \Alt (\cA)$ is not an annulus.
\hfill
$\Box$

\medskip
To show that if $\Omegalt$ is a topological annulus, then
$\la > 1/\ell$ we will need the following result.  

\begin{lemma}
\label{nonconstant-l}
Assume that $\la \leq 1/\ell$.
Let $\tau : \T \rightarrow \R$ be a continuous function.
If $\rho^+ - \rho^-$ is a constant function, then
$\Omegalt$ is not an annulus.
\end{lemma}

\proof
Suppose that  $\la \leq 1/\ell$ and $\rho^+ - \rho^-$ is the constant function  $ C$ for some $C >0$.
We proceed by contradiction.
If $\Omegalt$ is an annulus, then 
$$[\rho^- (\theta) , \rho^+(\theta)] = \cup_{\thetap \in \m^{-1}_\ell(\theta)} J(\thetap)$$
where $J(\thetap) = [\la \rho^- (\thetap) + \tau(\thetap),  \la \rho^+ (\thetap) + \tau(\thetap)].$
Since the length of $[\rho^- (\theta) , \rho^+(\theta)]$ is $C$ and the length of 
each of the $\ell$ intervals  $J(\thetap)$ is  $\la C$, it follows that $\ell \la C \geq C$ and therefore $\la = 1/\ell$.
Moreover, the interior of the intervals $J(\thetap)$ must be  pairwise disjoint, which is impossible
since the image under $\Alt$ of the loop $\{ (\theta, (\rho^+(\theta) + \rho^-(\theta))/2) \,\, / \, \theta \in \T \}$ 
must have self-intersections. 
\hfill $\Box$ 
          
\begin{proposition}
  \label{ell-p}
  If $\Omegalt$ is a closed topological annulus, then $\la > 1/\ell$.
\end{proposition}

\proof
Suppose that $\la \leq 1/\ell$.
Let $C$ be the maximum of $\rho^+ - \rho^-$ and let $E \subset \T$ be the set  formed by the arguments 
$\theta$ such that $C=(\rho^+ - \rho^-)(\theta)$.
>From the previous lemma we may assume that $E \neq \T$. 
Hence, there exists $\theta_0 \notin E$ such that $\theta = \ell \theta_0 \in E$.
It follows that the length of 
$$S= \cup_{\thetap \in \m^{-1}_\ell(\theta)} [\la \rho^- (\thetap) + \tau(\thetap),  \la \rho^+ (\thetap) + \tau(\thetap)]$$
is strictly less than $C$.
Therefore,  $[\rho^- (\theta) , \rho^+(\theta)]  \neq S$ and by Remark~\ref{annulus-r} we conclude
that $\Omegalt$ is not an annulus.
\hfill $\Box$

\subsection{Annular attractors when $\tau$ is not cohomologous to $0$}
\label{not-coboundary-ss}
The aim of this subsection is to prove (2) of Theorem~\ref{affine-th}
under some extra assumptions on $\tau$.

\begin{definition}
  \label{coboundary-d}
  We say that  $\tau: \T \rightarrow \R$ is {\bf cohomologous to $0$} or
a coboundary if $\tau = \mu \circ \m_\ell - \mu$ for some
continuous $\mu : \T \rightarrow \R$.
\end{definition}

We will show that Theorem~\ref{affine-th} (2) holds if $\tau$  is not cohomologous to $0$ and $\int_\T  \tau(\theta) d\theta =0$.
Under these  stronger assumptions we can weaken the Lipschitz class  hypothesis for $\tau$:

\begin{proposition}
\label{not-coboundary-p}
Let  $\tau: \T \rightarrow \R$ be a H\"older function which is not cohomologous to $0$ and such that 
$\int_\T \tau(\theta) d\theta = 0$.
 Then there exists $\lambdanot < 1$ such that $\Omegalt$ is a closed topological  annulus for all $\lainlazoo$.
\end{proposition}

The proof relies on finding appropriate periodic orbits under the multiplication by $\ell$ map 
$\m_\ell : \T \rightarrow \T$. 

\smallskip
To simplify notation we let $\int \tau = \int_\T  \tau(\theta) d\theta =0$ since we will only consider integrals 
with respect to the Lebesgue measure on $\T$.

\begin{lemma}
\label{subcohomological-l}
If $\tau: \T \rightarrow \R$ is a H\"older function 
which is not cohomologous to $0$ and such that $\int \tau = 0$,
then 
there exist periodic points $\theta_\pm$ of period $p_\pm$  such that:
\begin{eqnarray*}
  \tau(\theta_+) + \cdots + \tau(\ell^{p_+}  \theta_+) & > & 0 \\
  \tau(\theta_-) + \cdots + \tau(\ell^{p_-}  \theta_-) & < & 0.
\end{eqnarray*}
\end{lemma}

We employ the ideas contained in  the proof of Theorem 9 in~\cite{contreras-01}.

\smallskip
\proof
By considering $-\tau$ instead of $\tau$ it is sufficient to find $\theta_+$.
For $\btheta = (\theta_0, \theta_1, \theta_2, \dots ) \in \cS$ we let
$$S_n \btheta = \tau(\theta_1) + \cdots + \tau(\theta_n).$$
We proceed by contradiction and suppose that for all period $n$ points
$\btheta$ we have that $S_n \theta \leq 0$.

We claim that 
$$\{ S_n \btheta \, / \, \btheta \in \cS, n \in \N \}$$
is bounded above.
In fact, for any $\btheta \in \cS$ there exists
a period $n$ point $\bthetap = (\thetap_0, \thetap_1, \dots ) \in \cS$
such that $$\Tdist(\thetap_n , \theta_n ) \leq \frac{1}{\ell^{n-1}}.$$
Therefore,
$$S_n (\btheta) \leq 
S_n (\btheta) - S_n (\bthetap)  \leq C \cdot (1 + \cdots + \ell^{-(n-1)\alpha}) \leq C \cdot \frac{\ell^\alpha}{\ell^\alpha -1}$$
where $C >0$ is an  $\alpha$--H\"older constant constant for $\tau$. 

Now  $\mu : \T \rightarrow \R$  defined by
$$\mu (\theta_0) = \sup \{ S_n \btheta \, / \, \btheta = (\theta_0, \beta_1, \beta_2, \dots ) \in \cS \}$$
is an $\alpha$--H\"older function. 
In fact, for any $\theta_0, \thetap_0 \in \T$ and $\varepsilon >0$ let
$\btheta = (\theta_0, \theta_1, \dots)$ be such that $\mu (\theta_0) < S_n \btheta + \varepsilon$.
Then there exists $\bthetap = (\thetap_0, \dots) \in \cS$ such that 
$\Tdist(\theta_k, \thetap_k) = \ell^{-k} \Tdist(\theta_0, \thetap_0)$.
Therefore, 
\begin{eqnarray*}
  \mu(\thetap_0) \geq S_n (\bthetap) & = & S_n (\bthetap) - S_n (\btheta) + S_n (\btheta) \\
& \geq & S_n (\btheta) - C \frac{1}{\ell^\alpha-1} 
\Tdist (\theta_0, \thetap_0)^\alpha \\
& \geq & \mu(\theta_0) - \varepsilon -  C \frac{1}{\ell^\alpha-1} 
\Tdist (\theta_0, \thetap_0)^\alpha.
\end{eqnarray*}
It follows that $\mu$ is $\alpha$-H\"older, in particular continuous.

>From the definition of $\mu$ we conclude that
$$\mu (\ell \theta) \geq \mu (\theta) + \tau (\theta).$$
Therefore $\mu \circ \m_\ell - \mu - \tau$ is a non-negative function
whose integral over $\T$ is zero. It follows that
$\tau = \mu \circ \m_\ell - \mu$ which contradicts our assumption that
$\tau$ is not cohomologous to $0$.
\hfill
$\Box$ 

\medskip
\noindent
{\bf Proof of Proposition~\ref{not-coboundary-p}.}
For $\la$ sufficiently close to $1$ we will show that  $\rho^+ > \lambda^{-1} \|\tau\|_\infty$.
Similarly, one can  show that for $\la$ close to $1$ the inequality  $\rho^- < - \lambda^{-1} \|\tau\|_\infty$ holds.
In view of  Lemma~\ref{annulus-characterization-l}  this is sufficient to guarantee that
$\Omegalt$ is an annulus.

Let $\theta_+$ be a period $p=p_+$ periodic point as  in the previous lemma
and let $c = \tau(\theta_+) + \cdots + \tau(\ell^{p-1} \theta_+) >0$.  
Consider $\la_1 <1$ and a neighborhood $U \subset \T$ of $\cO = \{ \theta_+ , \dots, \ell^{p-1} \theta_+\} $ such that:

(a) There exists a neighborhood $V \subset U$ of $\cO$ such that $\m_\ell : V \rightarrow U$ is
a bijection.

(b) For all $\la > \la_1$, if $\{ \theta, \ell \theta, \dots, \ell^{p-1} \theta \} \subset U$, then 
$$\tau(\theta) + \la \tau (\ell \theta) + \cdots + \la^{p-1} \tau(\ell^{p-1} \theta) > \frac{c}{2}.$$

Let $n_0$ be such that $\m^{n_0}_\ell (U) = \T$.
Then for all $\theta_0 \in \T$ there exists $\btheta = (\theta_0 , \dots ) \in \cS$ such that
$\theta_n \in U$ for all $n \geq n_0$. 
Hence for all $\la > \la_1$,
\begin{eqnarray*}
  \rho^+ (\theta_0) & \geq & \sum_{k \geq 1} \la^{k-1} \tau (\theta_k)  =  \sum^{n_0-1}_{k \geq 1} \la^{k-1} \tau (\theta_k)
+ \sum_{k \geq n_0} \la^{k-1} \tau (\theta_k)\\
& \geq & -(1 + \cdots + \la^{n_0 -1}) \| \tau \|_\infty + \la^{n_0-1} (1 + \la^p + \la^{2p} + \cdots ) \frac{c}{2}\\
& \geq & - n_0 \| \tau \|_\infty  + \frac{c \la^{n_0-1}}{2(1-\la^p)}.
\end{eqnarray*}
It follows that there exists $\la_0 \in (\la_1 , 1)$
such that  for all $\lainlazoo$  we have that $\rho^+ > \lambda^{-1} \|\tau\|_\infty$.
\hfill
$\Box$ 

\begin{remark}{\em
\label{strict-r}
Provided that $\tau : \T \rightarrow \R$ is H\"older, $\int \tau =0$ and
$\tau$ is not cohomologous to zero, it follows from the previous proof that
for $\la \in [\la_0, 1)$ 
the inequality~(\ref{annular-e}) holds strictly, that is: 
\begin{equation*}
\min \{ \la \rho^+ (\thetap) + \tau (\thetap) \, / \, \thetap \in \m^{-1}_\ell (\theta) \} 
 >
\max \{ \la \rho^- (\thetap) - \tau (\thetap) \, / \, \thetap \in \m^{-1}_\ell (\theta) \}. 
\end{equation*}}
\end{remark}

\section{Cohomological operators}
\label{linear-s}
In this section we will study the linear operators
$$\cL_\la  \mu = \mu \circ \m_\ell - \la \mu$$
acting on Lipschitz functions $\mu : \T \rightarrow \R$
where $\lainzoc$.
On one hand these linear operators are related to 
conjugacy classes of affine maps $\Alt$ (see Lemma~\ref{affine-conjugacy-l}). On the
other hand, in view of 
propositions~\ref{jordan-characterization-p} and~\ref{not-coboundary-p}, these operators are also
 related to topological properties of
the attractor $\Omegalt$. At the end of this section we will apply the properties
of $\cL_\la$ to prove a stronger version of Theorem~\ref{affine-th} announced in Subsection~\ref{coho-ss}.

\medskip
The key to study certain topological properties of $\Omegalt$ is to understand
the iterated ``cohomological'' equations $\cL^n_\la \mu = \tau$. 

\begin{definition}
\label{multiplicity-d}
Consider a continuous function $\tau : \T \rightarrow \R$.
Given $\la \in (0,1]$ we define the multiplicity $\multt(\la)$ as follows:

$\bullet$ $\multt(\la) = 0$, if $\cL_\la \mu \neq \tau$ for all 
continuous $\mu$.

$\bullet$ $\multt(\la)= n$, if $\cL^n_\la \mu =  \tau$ for some integer $n >0$ and some 
continuous $\mu$  such that $\operatorname{mult}_{\mu - \int\mu} (\la) =0$.

$\bullet$ $\multt(\la)= \infty$, otherwise. 
\end{definition}

We will show that $\cJ^\prime_\tau = \{ \la \in (0,1] \,\,/\, \multt (\la) >0 \}$ is finite (counting multiplicities) provided that
$\tau$ is a non-constant Lipschitz function (cf. Lemma~\ref{finite-jordan-l}). In particular, $\multt(\la)$ is finite for each
$\la \in (0,1]$. 

\medskip
Although our main interest is on maps $\Alt$ where $\tau : \T \rightarrow \R$ is Lipschitz,
the basic properties of the operators $\cL_\la$ acting on any class are summarized in 
the following lemma.

\begin{lemma}
  \label{lla-l}
  For $\la \in (0 , 1]$, let  $\cL_\la : \cC \rightarrow \cC$ be the linear operator defined above acting on 
the class $\cC$ of maps from $\T$ to $\R$ where $\cC$ is the class of Lipschitz or $\cC^r$ maps 
for some $r \in [0, \infty] \cup \{ \omega \}$. Then:

{\em (1)} $\cL_{\la_1} \circ \cL_{\la_2} = \cL_{\la_2} \circ \cL_{\la_1}$ for all $0 < \la_1 , \la_2 \leq 1$.

{\em (2)} If $\cL_{\la_1} \mu_1 = \tau$, then 
$\operatorname{mult}_{\mu_1}(\la) = \multt(\la)$
for all $\la \neq \la_1$ and 
$\operatorname{mult}_{\mu_1} (\la_1) = \multt(\la_1) -1$ provided that $\int \mu_1 =0$ whenever $\la_1 =1$. 

{\em (3)} $\operatorname{ker}\cL_1 = \{ \mbox{constant maps} \}$ and 
$\operatorname{ker}\cL_\la = \{ 0 \}$ for all $\lainzoo$.

{\em (4)} If  $\cL_\la \mu = \tau$ and $\int \tau =0$ for some $\la \in (0,1)$, then
$\int \mu =0$.

{\em (5)} If $\tau : \T \rightarrow \R$ is of class $\cC$ and $\mu : \T \rightarrow \R$ is
a continuous function such that $\mu \circ \m_\ell - \la \mu = \tau$, then $\mu$ is of class $\cC$.

{\em (6)} If $\tau \in \operatorname{Image}\cL_1$, then $\int \tau  =0$.

{\em (7)} $\cL_\la (\{ \mbox{constant maps} \}) = \{ \mbox{constant maps} \}$ for all $\la \in (0,1)$.
\end{lemma}

\proof
Statement (1)  is a straightforward computation and statements (4) and (6) are an immediate consequence of the 
fact that the Lebesgue measure in the circle is $\m_\ell$--invariant (i.e., $\int \tau \circ \m_\ell = \int \tau$ for all
continuous $\tau : \T \rightarrow \R$). 

\smallskip
For (2), suppose that $\cL_{\la_1} \mu_1 = \tau$
and just note that if, for some $\la \neq \la_1$, there exists $\mu$ such that
 $\cL_{\la} \mu  = \tau$, then $\cL_{\la} (\mu - \mu_1)/(\la -\la_1)
=\mu_1$.  

\smallskip
(3) If $\mu (\ell \theta) -\la \mu(\theta) =0$, then $\mu (\ell^n \theta) = \la^n \mu (\theta)$ for all $\theta \in \T$
and all $n \geq 0$. Let $\theta_0$ be such that $\{ \ell^n \theta_0 \}_{n \geq 0}$ is dense in $\T$.
Hence, for all $\theta \in \T$,  $\mu (\theta) = 0$ when  $\la <1$ and $\mu (\theta) = \mu (\theta_0)$ when $\la =1$.

\smallskip
(5) As in the proof of Proposition~\ref{jordan-characterization-p} we pass to the universal
cover $\pi : \R \rightarrow \T$. That is, if  $\mu \circ \m_\ell - \la \mu = \tau$
then $\tmu (\ell s) - \la \tmu (s)= \ttau (s)$ where $\tmu = \mu \circ \pi$ and $\ttau = \tau \circ \pi$.
It follows that $\tmu (s) = \ttau (\ell^{-1} s) + \la \ttau (\ell^{-2} s) + \cdots$ is a class $\cC$ map
from $\R$ to $\R$ and therefore $\mu : \T \rightarrow \R$ is also of class $\cC$. 

\smallskip
Since the image under $\cL_\la$ of the constant function $\tau$ equal to  $1$ is the constant function
$\mu$ equal to $1 - \la$,  statement (7) follows.
\hfill $\Box$

\begin{remark}
\label{collection-r}
{\em 
>From (1) and (2) it follows that  
there exists a continuous $\mu$ such that
$\cL_{\la_1} \circ \cdots \circ \cL_{\la_m} \mu = \tau$
if and only if 
$\la_1, \dots, \la_m$ is a collection
of elements of $\cJ^\prime_\tau$, maybe with repetitions, but 
such that the number of occurrences of $\la \in \cJ^\prime_\tau$ is not greater than
its multiplicity.}
\end{remark}

Now we prove some basic conjugacy relations among affine maps of $\cyl$.

\begin{lemma}
  \label{affine-conjugacy-l}
  Let $\tau : \T \rightarrow \R$ be a map of class $\cC$ where $\cC$ is the class 
of Lipschitz or $\cC^r$ maps for some  $r \in (0, \infty] \cup \{ \omega \}$. Then:

{\em (1)} For all $c \in \R$ and $0 < \la <1$, the maps $\Alt$ and $A_{\la, \tau +c }$ are
conjugate via an affine map.

{\em (2)} If $\mu \circ  \m_\ell - \la_0 \mu = \tau$ for some continuous map $\mu : \T \rightarrow  \R$
and $\la_0 \in (0 ,  1]$, then $\Alt$ and $A_{\la, \mu}$ are $\cC$--conjugate for all $\la \neq \la_0$. 
\end{lemma}

\proof
For (1) note that $\Alt \circ h = h \circ A_{\la, \tau +c}$ where  $h \tht = (\theta, t - c/(1-\la))$.

For (2), if $\mu \circ  \m_\ell - \la_0 \mu = \tau$ for some continuous map $\mu$, then $\mu$ is automatically
of class $\cC$ and $h \tht = (\theta, (\la - \la_0) t + \mu (\theta))$ is such that
$h \circ A_{\la, \mu} = \Alt \circ h$ for all $\la \neq \la_0$.
\hfill
$\Box$

\begin{lemma}
\label{lip-norm-l}
If $\mu : \T \rightarrow \R$ is Lipschitz, then
$$||\cL_\lambda \mu ||_L  \geq (\ell - \lambda) || \mu ||_L$$
for all $\la \in (0,1]$. 
\end{lemma}
 
\proof
Let $\mu$ be Lipschitz and let $\tau = \cL_\la \mu$.
Consider $\btheta = (\theta_k) \in \cS$ and observe that
$$\mu (\theta_0)= \tau(\theta_1) + \la \tau(\theta_2) + \cdots + \la^n\tau(\theta_{n+1}) - \la^{n+1} \mu (\theta_{n+1}).$$
For $\thetap_0 \in \T$ let $\bthetap = (\thetap_k)$ be such that
$$\Tdist (\theta_k , \thetap_k) = \frac{\Tdist(\theta_0, \thetap_0)}{\ell^k}.$$
It follows that for all $n \geq 0$:
\begin{eqnarray*}
  |\mu(\theta_0) - \mu(\thetap_0)| & \leq & \Tdist(\theta_0, \thetap_0) || \tau ||_L (\frac{1}{\ell} + \frac{\la}{\ell^2} + \dots
+ \frac{\la^n}{\ell^{n+1}}) + \\
& & \mbox{ } \la^{n+1} |\mu (\theta_{n+1}) - \mu(\thetap_{n+1})| \\
                                   & \leq & \Tdist(\theta_0, \thetap_0) || \tau ||_L \frac{1}{\ell - \la} + 
                                                                          \la^{n+1} |\mu (\theta_{n+1}) - \mu(\thetap_{n+1})|. 
\end{eqnarray*}
Therefore, $$ |\mu(\theta_0) - \mu(\thetap_0)| \leq || \tau ||_L \frac{1}{\ell - \la}  \Tdist(\theta_0, \thetap_0)$$
and the lemma follows.
\hfill
$\Box$

\begin{lemma}
\label{closed-lla-l}
Let $\tau : \T \longrightarrow \R$ be a Lipschitz map such that $\int \tau =0$.
The set 
$$\cJ^\prime_\tau = \{ \la \in (0,1] \, / \, \multt(\la) > 0 \}$$
is closed.
\end{lemma}

Note that from Proposition~\ref{jordan-characterization-p} it follows that
$$\cJ^\prime_\tau \cap (0,1) = \cJ_\tau = \{ \la \in (0,1) \, / \, \mbox{$\Omegalt$ is a Jordan curve} \}.$$

\medskip
\proof
Suppose that there exists $\la_n \rightarrow \la \in (0,1]$ and $\mu_n : \T \rightarrow \R$
such that $\mu_n \circ \m_\ell - \la_n \mu = \tau$. 
It follows that $\int  \mu_n =0$ and $\| \mu_n \|_L \leq \| \tau \|_L$
for all $n$. Therefore, $\{ \mu_n \}$ is an equicontinuous and uniformly bounded family.
Hence, by passing to a subsequence, we may assume that $\{ \mu_n \}$ converges 
to some continuous function $\mu$ which necessarily satisfies the equation $\mu \circ \m_\ell - \la \mu = \tau$,
that is,  $\la \in \cJ^\prime_\tau$.
\hfill $\Box$

\medskip
We now show that a Lipschitz function $\tau$ which is not constant 
is not ``infinitely'' cohomologous to $0$ or equivalently that $\la = 1$ has finite multiplicity
(see Definition~\ref{multiplicity-d}).

\begin{lemma}
\label{finite-boundary-l}
Let $\tau : \T \longrightarrow \R$ be a non-constant Lipschitz function such that $\int \tau =0$. 
Then there exists $m \geq 0$ and a Lipschitz function $\mu:\T \rightarrow \R$ such that
$\cL^m_1 \mu = \tau$ and $\int \mu =0$ but 
 $\cL_1 \rho  \neq  \mu$ for all Lipschitz functions $\rho$. 
\end{lemma}

\proof
Suppose that for $0 \leq n \leq m$ there exist  $\mu_n : \T \rightarrow \R$ 
such that $$\cL^n \mu_n = \mu_0 = \tau.$$
By Lemma~\ref{lla-l}, maybe after adding a constant to $\mu_m$,
we may assume that $\int \mu_n =0$ for $0 \leq n \leq m$. It follows that  $\cL_1 \mu_n = \mu_{n-1}$.
Under the assumption that $\tau = \mu_0$ is not identically $0$ we will exhibit an upper bound for 
$m$ in terms of the Fourier coefficients of $\tau$.

The Fourier coefficients
$\hmu_n(k)$ are uniformly bounded.
In fact, by Lemma~\ref{lip-norm-l}, $\| \mu_n \|_L \leq \| \tau \|_L$ and therefore 
\begin{equation}
  | \hmu_n (k) | \leq \frac{\| \tau \|_L}{4k}
\label{fourier-e}
\end{equation}
 for all $0 \neq k \in \Z$ and  $0 \leq n \leq m$.

Since $\mu_0$ is not identically $0$, there exists $k \in \Z$ such that $\ell \nmid k$ and $p \geq 0$ such that:
\begin{eqnarray}
  \hmu_0 (\ell^j k) & = & 0 \,\,\,\,\, \mbox{ for } 0 \leq j < p, \nonumber \\
  \hmu_0 (\ell^p k) & \neq & 0. \nonumber
\end{eqnarray}

Taking the $\ell^j k$--th Fourier coefficient to $\mu_n \circ \m_\ell - \mu_n = \mu_{n-1}$:
\begin{eqnarray}
\hmu_n ( \ell^{j-1} k) -  \hmu_n (\ell^j k)& = & \hmu_{n-1} (\ell^j k) \mbox{ for } j \geq 1,   \label{fourier-1-ej}\\
- \hmu_n (k) & = & \hmu_{n-1} (k),   \label{fourier-1-e0}
\end{eqnarray}
for $1 \leq n \leq m$. 

By induction in $p \geq 0$, it is easy to deduce from  (\ref{fourier-1-ej})  and (\ref{fourier-1-e0}) 
that
if $\hmu_0 (\ell^j k) =0$ for $0 \leq j <p$, then:
\begin{eqnarray*}
\hmu_n (\ell^j k) &=&0 \mbox{ for } 0 \leq j < p, \\
-\hmu_n (\ell^p k) & = &\hmu_{n-1}(\ell^p k),
\end{eqnarray*}
 for $0 \leq n \leq m$.
Therefore, from~(\ref{fourier-1-ej}), it follows that
$$m (-1)^{n-1} \hmu_0 (\ell^p k) - \hmu_m(\ell^{p+1} k) = (-1)^{n+1} \hmu_0 (\ell^{p+1} k).$$
Hence,
$$m \leq \frac{|\hmu_0 (\ell^{p+1} k)| + |\hmu_m(\ell^{p+1} k)|}{ |\hmu_0 (\ell^p k)|}$$
and by~(\ref{fourier-e}) we obtain an upper bound for $m$:
\begin{equation}
m \leq \frac{ \|\tau\|_L + 4k\ell^{p+1}|\htau (\ell^{p+1} k)|}{4k\ell^{p+1} |\htau (\ell^p k)|}.
\label{mbound-e}
\end{equation}
\hfill
$\Box$ 

\medskip
Below we record the explicit bound obtained in the previous proof.

\begin{corollary}
  Let $\tau: \T \rightarrow \R$ be a Lipschitz function such  that:
\begin{eqnarray}
  \htau (\ell^j k) & = & 0 \,\,\,\,\, \mbox{ for } 0 \leq j < p \nonumber \\
  \htau (\ell^p k) & \neq & 0 \nonumber
\end{eqnarray}
  for some integers $k \in \Z$  and $p \geq 0$ with $\ell \nmid k$.
If there exists $m \geq 0$ and a Lipschitz function $\mu : \T \rightarrow \R$ such that
$$\cL^m_1 \mu = \tau,$$
Then 
\begin{equation}
m \leq \frac{ \|\tau\|_L + 4k\ell^{p+1}|\htau (\ell^{p+1} k)|}{4k\ell^{p+1} |\htau (\ell^p k)|}. \nonumber
\end{equation}
\end{corollary}

\begin{lemma}[Main]
\label{finite-jordan-l}
Let $\tau : \T \longrightarrow \R$ be a Lipschitz function 
which is not constant.
Then $\cJ^\prime_\tau$ is finite, counting multiplicities.

Moreover, consider 
the finite collection $0 < \lambda_1, \ldots, \lambda_m \leq 1$ 
consisting of elements of  $\cJ^\prime_\tau$ where the number
of repetitions of each element of $\cJ^\prime_\tau$ coincides with its multiplicity.
Then there exists a Lipschitz function $\mu : \T \rightarrow \R$ 
such that
$$
\tau = \cL_{\lambda_1} \circ \ldots \circ \cL_{\lambda_m} \mu,
$$
and $\mu \neq \cL_\lambda \rho$  for all Lipschitz maps
$\rho : \T \rightarrow \R$ and all $\la \leq 1$. Furthermore, if $\int \tau = 0$, then
$\mu$ may be chosen so that $\int \mu =0$.
\end{lemma}

\proof
By Lemma~\ref{lla-l} and Remark~\ref{collection-r} it is sufficient to show that $\cJ^\prime_\tau$ is finite, counting multiplicities.
In view of  Lemma~\ref{lla-l} (7),  $\multt(\la) = \operatorname{mult}_{\tau -c} (\la)$ for all $0 < \la < 1$ and $c \in \R$, 
and since the multiplicity of $\la =1$
is always finite (cf. Lemma~\ref{finite-boundary-l}) it follows that 
 $\cJ^\prime_\tau$ is finite if and only if $\cJ^\prime_{\tau-\int \tau}$ is.
Therefore we may assume that $\int \tau =0$.

We proceed by contradiction and suppose 
$\cJ^\prime_\tau$ is infinite or that it contains an element
with infinite multiplicity. In both cases,  for all $n \geq 1$ there exist 
$0< \la_n \leq 1$ and $\mu_n : \T \rightarrow \R$ Lipschitz such that  
\begin{equation}
\cL_{\la_1} \circ \cdots \circ \cL_{\la_n} \mu_n = \tau.
\label{comp-e}
\end{equation}
By the previous lemma and by Lemma~\ref{lla-l} (1), 
the number of occurrences of 
$\la_n =1$ is finite. 
In particular, there exists $N \geq 0$ such that $0<\la_n < 1$ for all $n >N$
and $$\cL_{\la_{N+1}} \circ \cdots \circ \cL_{\la_{N+k}} (\cL_{\lambda_1} \circ \ldots \circ \cL_{\lambda_N} \mu_n) = \tau. $$
Therefore, after relabelling, we may assume that for all $n \geq 1$ there exist 
$\la_n < 1$ and $\mu_n : \T \rightarrow \R$ Lipschitz such that (\ref{comp-e}) holds.
Note that since $\int \tau =0$, then $\int \mu_n =0$ for all $n \geq 1$ (see Lemma~\ref{lla-l} (4)).

By Lemma~\ref{closed-lla-l}, it follows that there exists $\la_0<1$ such that
$\la_n \leq \la_0$ for all $n \geq 1$. 
Since $$||\tau||_L \geq (\ell - \la_1) \cdots (\ell - \la_n) ||\mu_n ||_L,$$
we have that $||\mu_n||_L \leq ||\tau||_L$ for all $n$. Therefore,
for all $0 \neq k \in \Z$ and $n \geq 1$,
$$|\hmu_n (k) | \leq ||\tau||_L,$$
where $\hmu_n(k)$ denotes the $k$-th Fourier coefficient of $\mu_n$.
Now by Lemma~\ref{lla-l} (3), for all $n \geq 1$
$$\cL_{\la_n} \mu_n = \mu_{n-1}$$ 
where $\mu_0 = \tau$.
Hence, $\mu_n (\ell \theta) - \la_n \mu_n (\theta) = \mu_{n-1}(\theta)$.
Therefore, for all $0 \neq k \in \Z$ such that $\ell \nmid k$:
\begin{eqnarray}
\hmu_n ( \ell^{j-1} k) - \la_n \hmu_n (\ell^j k)& = & \hmu_{n-1} (\ell^j k) \mbox{ for } j \geq 1,   \label{fourier-la-ej}\\
-\la _n \hmu_n (k) & = & \hmu_{n-1} (k).   \label{fourier-la-e0}
\end{eqnarray}

Fix $k \neq 0$ such that $\ell$ does not divide $k$.
We show by induction on $j \geq 0$ that $\hmu_n ( \ell^j k) =0$ for all
$n \geq 0$.

For $j=0$, by equation (\ref{fourier-la-e0}), $\hmu_0 (k) = (-1)^n \la_1 \cdots \la_n \hmu_n (k)$.
Since $\hmu_n(k)$ are uniformly bounded it follows that $\hmu_0 (k) =0$ and therefore
$\hmu_n(k) =0$ for all $n$.

Let $j \geq 1$ and suppose that $\hmu_n ( \ell^{j-1} k) =0$ for all
$n \geq 0$. It follows from equation  (\ref{fourier-la-ej}) that $- \la_n \hmu_n (\ell^j k)  =  \hmu_{n-1} (\ell^j k)$.
Hence,  $\hmu_0 (\ell^j k) = (-1)^n \la_1 \cdots \la_n \hmu_n (\ell^j k)$. Again using that $\hmu_n (\ell^j k)$ are
uniformly bounded, it follows that $\hmu_n (\ell^j k)=0$ for all $n$. 

We conclude that $\tau = \mu_0$ is constant which is a contradiction.
\hfill $\Box$

\medskip
We now  state and prove a stronger version of Theorem~\ref{affine-th}.
Recall that  $\cJ_\tau$ denotes the set of $\la \in (0,1)$ such that
$\Alt$ is a Jordan curve. 

\begin{theorem}
\label{strong-affine-th}
Suppose that $\tau : \T \rightarrow \R$ is Lipschitz. 
Then there exists a Lipschitz map 
$\mu: \T \rightarrow \R$ satisfying the following properties:

{\em(1)} $\cJ_\mu = \emptyset$ and  for all $\la \in (0,1) \setminus \cJ_\tau$
the maps $\Alt$ and $A_{\la, \mu}$ are topologically conjugate.

{\em (2)}
$\int  \mu   =0 $ and $\mu$ is not cohomologous to $0$.
\end{theorem}

\medskip
\noindent
{\bf Proof of Theorem~\ref{strong-affine-th}.}
By Lemma~{\ref{affine-conjugacy-l}} (1), after replacing $\tau$ by $\tau - \int \tau$, we may assume that
$\int \tau =0$. In view of Proposition~\ref{jordan-characterization-p} and Lemma~\ref{affine-conjugacy-l} (2),
statements (1) and (2) hold for the function $\mu : \T \rightarrow \R$ given by the previous lemma.
\hfill
$\Box$

\medskip
\noindent
{\bf Proof of Theorem~\ref{affine-th}.}
For (1) just note that Lemma~\ref{finite-jordan-l} implies that $\cJ_\tau$ is finite.
For (2) let $\mu : \T \rightarrow \R$ be as in Theorem~\ref{strong-affine-th} and apply Proposition~\ref{not-coboundary-p}.
\hfill
$\Box$

\section{Continuity of the upper and lower boundaries.}
\label{uplow-s}
Perturbation of affine maps $\Alt$ also posses an attractor for which the upper and lower boundaries may be defined
(see Definition~\ref{uplowF-d} below).
The aim of this section is to show that the upper and lower boundaries of 
$\Omegalt$ vary continuously under Lipschitz perturbations of $\Alt$. 

Given an open set $U \subset \cyl$ with compact closure we denote by 
$\lip ( \overline{U}, \cyl)$ the set formed by all the Lipschitz maps
$F : \overline{U} \rightarrow \cyl$ endowed with the Lipschitz metric
$\distL$. More precisely, 
$$\distL (F_0, F_1) = \| F_0 - F_1 \|_\infty + \| f_0 -f_1\|_L + \|g_0 -g_1\|_L$$
where $F_i \tht = (f_i \tht , g_i \tht)$ for $i=0,1$. 

Note that the set $\lip(\overline{U},U)$  of all $F \in \lip ( \overline{U}, \cyl)$ such that
$F(\overline{U}) \subset U$ is open in  $\lip ( \overline{U}, \cyl)$.

\begin{definition}
\label{uplowF-d}
Suppose that $U \subset \cyl$ is an  open set with compact closure.
For any $F \in \lip (\overline{U}, U)$ we let $\Omega_F = \cap_{n \geq 0} F^n (\overline{U})$ be
the attractor of $F$. The {\bf upper and lower boundaries of $\Omega_F$} are the graphs of 
$$\rho^\pm_F : \T \rightarrow \R \cup \{\pm \infty\}$$
where
$$\rho^+_F (\theta) = \sup \{ t \in \R \,\, /\, \tht \in \Omega_F \}$$
$$\rho^-_F (\theta) = \inf \{ t \in \R \,\, /\, \tht \in \Omega_F \}$$
if $\Omega_F \cap (\{ \theta \} \times \R) \neq \emptyset$ 
and $\rho^\pm(\theta) = \mp \infty$ otherwise.
\end{definition}

We may now state the main result of this section.

\begin{proposition}
\label{uplow-p}
Consider $\la \in (0,1)$ and a Lipschitz function $\tau : \T \rightarrow \R$.
Let $T_0 \in \R$ be such that $\Alt(\overline{U_0}) \subset U_0$ where $U_0 = \T \times (-T_0,T_0)$.
Given $\epsilon >0$ there exists a neighborhood $\cU \subset \lip(\overline{U_0},U_0)$ of 
$\Alt$ such that for all $F \in \cU$ the following hold:

{\em (1)} $\rho^\pm_F : \T \rightarrow \R$ are well defined 
Lipschitz functions.

{\em (2)} $\| \rho^\pm_F - \rho^\pm_{\Alt} \|_\infty  < \epsilon.$
\end{proposition}

\subsection{Lipschitz maps in $\R^2$}
\label{r2-ss}
In this subsection we consider $\bV = \R \times [-T_0,T_0]
\subset \R^2$ and study the action of Lipschitz maps $\tF : \bV \rightarrow \R^2$ 
on the graphs of Lipschitz functions $\trho : \R \rightarrow [-T_0,T_0]$.  
The results obtained here will be applied in the proof of 
Proposition~\ref{uplow-p} to the lift of 
maps $F$ which are Lipschitz close to $\Alt$.

Throughout this subsection we consider a map:
$$\begin{array}{rccc}
\tF :& \bV & \rightarrow & \R^2 \\
& (s,t) & \mapsto &  ( \tf(s,t) , \tg(s,t))
\end{array}$$
and suppose that there exist positive constants $\ell_0, \la_0, C_{1,2}$ and $C_{2,1}$ such that:

\begin{eqnarray}
\tf (s_0,t) - \tf (s_1,t) & \geq & \ell_0 (s_0 - s_1), \label{cont1-e} \\
|\tg (s_0,t) - \tg   (s_1,t)| & \leq & C_{2,1} |s_0 - s_1|,  \label{cont3-e}\\
|\tf (s,t_0) - \tf   (s,t_1)| & \leq & C_{1,2} |t_0 - t_1|,  \label{cont4-e}\\
|\tg (s,t_0) - \tg   (s,t_1)| & \leq & \la_0 |t_0 - t_1|,  \label{cont5-e}
\end{eqnarray}
for all $s_0 > s_1$ and $t,t_0,t_1 \in [-T_0,T_0]$.

\smallskip
The above conditions are satisfied by the lift $\tAlt (s,t) = (\ell s, 
\la t + \tau \circ \pi (s))$ of $\Alt$ where $\tau : \T \rightarrow
\R$ is Lipschitz and $\la \in (0,1)$. Observe that in this case we may choose
$\ell_0 =  \ell$, $C_{2,1} = \| \tau \|_L$, $C_{1,2} =0$, and
$\la_0 = \la$. 

\begin{lemma}
\label{action-l}
If $\tgamma \subset \bV$ is the graph of a Lipschitz function 
$\trho : \R \rightarrow [-T_0,T_0]$ 
such that $\|\trho\|_L < \ell_0 C^{-1}_{1,2}$, then $\tF(\tgamma)$ is the graph 
of a $C$-Lipschitz function $\cT_\tF (\trho): \R \rightarrow \R$, where 
$$C = \frac{C_{2,1} + \la_0 \| \trho\|_L}{\ell_0 - C_{1,2} \| \trho \|_L}.$$
\end{lemma}

\proof
Suppose that $s_0 > s_1$. For  $i=0,1$,  let $t_i = \trho(s_i)$ and $(s^\prime_i,t^\prime_i)
=\tF (s_i,t_i)$. Then we have that:
\begin{eqnarray*}
s^\prime_0 - s^\prime_1 & = & \tf (s_0,t_0) - \tf(s_1,t_1) \\
                        & = & (\tf (s_0,t_0) - \tf(s_0,t_1)) + (\tf (s_0,t_1) - \tf(s_1,t_1)) \\
                        & \geq & -C_{1,2} | t_0 - t_1 | + \ell_0 (s_0 - s_1) \\
                        & > &- \ell_0 |s_0 - s_1| + \ell_0 (s_0 - s_1) = 0.
\end{eqnarray*}
In particular, $\tF(\tgamma)$ is the graph of some function.
Also,
\begin{eqnarray*}
|t^\prime_0 - t^\prime_1| & = & |\tg (s_0,t_0) - \tg(s_1,t_1)| \\
                          & = & |\tg (s_0,t_0) - \tg(s_0,t_1) + \tg (s_0,t_1) - \tg(s_1,t_1)| \\
                          & \leq & \la_0 |t_0 - t_1| + C_{2,1} | s_0 - s_1 |.
\end{eqnarray*}
Hence,
\begin{eqnarray*}
\frac{|t^\prime_0 - t^\prime_1|}{|s^\prime_0 - s^\prime_1|} 
& \leq & \frac{\la_0 |t_0 - t_1| + C_{2,1} | s_0 - s1 |}{-C_{1,2} | t_0 - t_1 | + \ell_0 |s_0 - s_1|} \\ 
&\leq & \frac{C_{2,1} + \la_0 \| \trho\|_L}{\ell_0 - C_{1,2} \| \trho \|_L}
\end{eqnarray*}
and the Lemma follows.
\hfill $\Box$

\begin{definition}
\label{preserves-d}
  We say that  a Lipschitz map $\tF : \overline{V} \rightarrow \R^2$ {\bf preserves $C$-Lipschitz graphs} 
with constants $\ell_0, \la_0, C_{1,2}$ and $C_{2,1}$ if  (\ref{cont1-e})--(\ref{cont5-e}) hold
and 
$$ \frac{C_{2,1} + \la_0 C}{\ell_0 - C_{1,2} C} \leq C < \frac{\ell_0}{C_{1,2}}.$$
\end{definition}

In particular, if $\tF : \overline{V} \rightarrow \overline{V} \subset \R^2$ preserves
$C$--Lipschitz graphs, then $\cT_\tF$ acts on the set of $C$--Lipschitz functions 
$\rho : \R \rightarrow [-T_0,T_0]$ (see  Lemma~\ref{action-l}).

\smallskip
Now we compute a Lipschitz constant for $\cT_\tF$ with respect to the $\cC^0$--norm.

\begin{lemma}
\label{tf-lip-l}
Let $C >0$ be such that $ C < \ell_0 C^{-1}_{1,2}$.
Suppose that  $\trho_i : \R \rightarrow [-T_0,T_0]$ 
are $C$-Lipschitz functions where $i=0,1$.
Then
$$\| \cT_\tF (\trho_0) - \cT_\tF (\trho_1) \|_\infty \leq \frac{\la_0 \ell_0 + C_{2,1} C_{1,2}}{\ell_0 - C_{1,2} C}
\|\trho_0 - \trho_1 \|_\infty.$$
\end{lemma}

\proof
Consider $s^\prime \in \R$.
For $i=0,1$, let $\trho^\prime_i = \cT_\tF (\trho_i)$,
 $t^\prime_i = \trho^\prime_i (s^\prime)$.
Also we let  $(s_i,t_i)$ be the points in the graph of $\trho_i$ such that
$\tF(s_i,t_i) = (s^\prime,t^\prime_i)$. We may assume that $s_1 > s_0$.

We must find an upper bound for $|t^\prime_1 - t^\prime_0 |$ in terms of
$|\hat{t}_1 - t_0 |$ where $\widehat{t}_1 = \trho_1 (s_0)$.

Since 
$$0=  \tf (s_0,t_0) - \tf(s_1,t_1) =  \tf (s_0,t_0) - \tf(s_1,t_0) +  \tf (s_1,t_0) - \tf(s_1,t_1),$$
we have that
$$ \tf (s_1,t_0) - \tf(s_0,t_0) =  |\tf (s_1,t_1) - \tf(s_1,t_0)|.$$
Therefore,
\begin{equation}
\ell_0 (s_1 - s_0) \leq C_{1,2} |t_1 - t_0|. \label{tf1-e}
\end{equation}

Also,
\begin{eqnarray}
|t^\prime_1 - t^\prime_0 | & = & |  \tg (s_0,t_0) - \tg(s_1,t_1)| \nonumber \\
                           & \leq  & | \tg (s_0,t_0) - \tg(s_0,t_1)| + | \tg (s_0,t_1) - \tg(s_1,t_1)| \nonumber \\
                           & \leq & \la_0  |t_0 - t_1| + C_{2,1} | s_0 - s_1 |. \label{tf2-e} 
\end{eqnarray}

Since
\begin{eqnarray*}
|t_1 - t_0| & \leq &  |\hat{t}_1 - t_0 | + |\hat{t}_1 - t_1| \\
            & \leq &  |\hat{t}_1 - t_0 | + C | s_0 - s_1 |,
\end{eqnarray*}
it follows from (\ref{tf1-e}) that
$$
|t_1 - t_0 | \leq  |\hat{t}_1 - t_0 | + C C_{1,2} \ell^{-1}_0 |t_1 - t_0 |.
$$                                                                                                             
Hence,
\begin{equation}
|t_1 - t_0 | \leq (1 - C C_{1,2} \ell^{-1}_0)^{-1} |\hat{t}_1 - t_0 |. \label{tf5-e}
\end{equation}
Combining (\ref{tf1-e}) and (\ref{tf2-e}) we obtain
\begin{equation}
|t^\prime_1 - t^\prime_0 | \leq (\la_0 + C_{2,1}C_{1,2} \ell^{-1}_0) |t_1 - t_0 |. \label{tf4-e}
\end{equation}
The lemma now follows directly from (\ref{tf5-e}) and (\ref{tf4-e}).
\hfill
$\Box$

\begin{lemma}
\label{tfie-l}
For $i=0,1$, 
let $\tF_i(s,t) = (\tf_i (s,t), \tg_i(s,t))$ be maps in $\lip(\overline{V}, \R^2)$ 
such that $\tF_i$ preserve $C$-Lipschitz graphs 
with constants $\ell_0, \la_0, C_{1,2}$ and $C_{2,1}$.
If $\|\tF_0 - \tF_1 \|_\infty < \epsilon$, then 
$$\|\cT_{\tF_0} (\trho) - \cT_{\tF_1}(\trho)\|_\infty \leq (1 + C) \epsilon$$
for all $C$-Lipschitz functions $\trho: \R \rightarrow [-T_0,T_0]$. 
\end{lemma}

\proof
Consider $s_0 \in \R$ and let $t_0 = \trho(s_0)$.
Also let $(s^\prime_i, t^\prime_i) = F_i (s_0, t_0)$ where $i=0,1$.
It follows that
\begin{eqnarray*}
| \trho^\prime_1 (s^\prime_1) - \trho^\prime_0 (s^\prime_1)| & \leq & | \trho^\prime_1 (s^\prime_1) - \trho^\prime_0 (s^\prime_0)| + |\trho^\prime_0 (s^\prime_0) - \trho^\prime_0 (s^\prime_1)| \\
& \leq & \epsilon + C | s^\prime_0 - s^\prime_1| \leq \epsilon + C \epsilon
\end{eqnarray*}
where $\cT_{F_i}(\trho)=\trho^\prime_i$.
\hfill
$\Box$

\subsection{The upper and lower boundaries as fixed points}
\label{fixed-ss}
Throughout this subsection we let $T_0 >0$ and $U_0 = \T \times (-T_0, T_0)$.
Here we consider a map $F \in \lip ( \overline{U_0},U_0)$
and show that under certain conditions the upper and lower boundaries of the attractor $\Omega_F$ are
Lipschitz graphs.

\begin{lemma}
Let $F \in \lip ( \overline{U_0},U_0)$ be locally an orientation preserving
homeomorphism which acts as multiplication by $\ell \geq 1$ on the first homology of $\T \times (-T_0, T_0)$.
Denote by $\tF : \R \times [-T_0, T_0] \rightarrow \T \times (-T_0, T_0)$ a lift of $F$ to the universal cover
and suppose that $\tF$ preserves $C$-Lipschitz maps
with constants $\ell_0$, $\la_0$, $C_{1,2}$, $C_{2,1}$.
Assume that:
$$ \frac{\la_0 \ell_0 + C_{2,1} C_{1,2}}{\ell_0 - C_{1,2} C} < 1.$$
Let $$\cT^+_F (\trho) (s) = \max \{ \cT_\tF (\trho)(s), \cT_\tF (\trho)(s+1), \dots ,  \cT_\tF (\trho)(s+\ell-1) \}$$
and
   $$\cT^-_F (\trho) (s) = \min \{ \cT_\tF (\trho)(s), \cT_\tF (\trho)(s+1), \dots ,  \cT_\tF (\trho)(s+\ell-1) \}.$$
Then $\cT^\pm_F$ are contractions in the space of $C$-Lipschitz functions  $\trho : \R \rightarrow [-T_0,T_0]$ endowed
with $\| \cdot \|_\infty$. The fixed points $\trho^\pm_F$ 
of  $\cT^\pm_F$ are $1$-periodic and the graphs of
$$\begin{array}{rccc}
\rho^\pm_F :& \T & \rightarrow & \R \\
& \theta = \pi(s) & \mapsto &  \trho^\pm_F (s)
\end{array}$$
are the upper and lower boundaries of the attractor $\Omega_F$.
\end{lemma}

\proof
Since the maximum and minimum of $C$-Lipschitz functions are also $C$-Lipschitz, from Lemma~\ref{action-l} it follows
that $\cT^\pm_F (\trho)$ are  $C$-Lipschitz whenever $\trho$ is $C$-Lipschitz.
By Lemma~\ref{tf-lip-l}, the operators $\cT^\pm_F$ act as contraction maps.
Our hypothesis that $F$ acts as multiplication by $\ell$ on the first homology group 
translates to the universal cover as $\tF (s+1,t) = \tF(s,t) + (\ell, 0)$.
It follows that $\cT^\pm_F$ preserve the closed subset of $1$-periodic $C$-Lipschitz function.
Therefore, the fixed points $\trho^\pm_F$ are $1$-periodic and we let $\rho^\pm_F (\pi (s)) = \trho^\pm_F (s)$.

We now show that $\rho^+_F$ is the upper boundary of $\Omega_F$.
Note that the graph $\gamma^+_F$ of $\rho^+_F$ is invariant under $F$.
That is $F(\gamma^+_F) \supset \gamma^+_F$. Therefore $\gamma^+_F \subset \Omega_F$.
It is sufficient to show that $\Omega_F$ is below $\gamma^+_F$.
For this let $\trho_0$ be the constant function  $T_0$ on $\R$ and let $\trho_n = \cT^+_F (\trho_0)$.
>From our previous discussion and the fact that  $\trho_0$ is $1$-periodic, we conclude that
 $\trho_n$ projects to a function $\rho_n: \T \rightarrow \R$ with graph $\gamma_n$. 
Since $F$ is locally an orientation preserving homeomorphism, 
$F^n (\overline{U_0})$ has as upper boundary the curve $\gamma_n$.
Taking into consideration that $\rho_n$ converges to $\rho$,
it follows that  $\cap F^n (\overline{U_0})$ is below $\gamma^+_F$.
A similar argument shows that $\rho^-_F$ is the lower boundary of $\Omega_F$.
\hfill
$\Box$

\medskip
The analogue of Lemma~\ref{annulus-characterization-l} also holds in this context. More precisely:

\begin{lemma}
\label{annularF-l}
Assume that  $F \in \lip ( \overline{U_0},U_0)$ is  locally an orientation preserving
homeomorphism such that the upper and lower boundaries of $\Omega_F$ are 
the graphs $\gamma^\pm_F$  of functions $\rho^\pm_F : \T \rightarrow \R$. If for all $\theta \in \T$
\begin{equation}
\min \{ t \,\, / \, \tht \in F(\gamma^+_F)  \} 
 >
\max \{  t \,\, / \, \tht \in F(\gamma^-_F) \},
\end{equation}
then $\Omega_F$ is a closed topological annulus.
Moreover, if $\ell =2$ and  for some $\theta \in \T$
\begin{equation}
\min \{ t \,\, / \, \tht \in F(\gamma^+_F)  \} 
 <
\max \{  t \,\, / \, \tht \in F(\gamma^-_F) \},
\end{equation}
then $\Omega_F$ is not a closed topological annulus.
\end{lemma}

We omit the proof of this Lemma since it is very similar to that of Lemma~\ref{annulus-characterization-l}.

\subsection{Proof of Proposition~\ref{uplow-p}.}
Since Lipschitz perturbations of bilipschitz maps are bilipchitz,
all $F$ in a  sufficiently small neighborhood $\cU$ of $\Alt$ in $\lip(\overline{U_0},U_0)$ 
are locally orientation preserving homeomorphisms.
Consider a small real number $\delta >0$ and shrink 
$\cU$, if necessary, so that for all $F \in \cU$ there exists a unique  lift
$\tF: \R \times [-T_0,T_0] \rightarrow \R^2$ of $F$  such  that
$\| \tF - \tAlt \|_\infty < \delta$, where $\tAlt (s,t) = (\ell s, \la t + \tau \circ \pi (s))$.
After further shrinking of $\cU$, if necessary, we may assume that there exist
positive constants: $\ell_0$ slightly smaller than $\ell$, $\la_0$ close to $\la$, 
$C_{2,1}$ close to $\|\tau\|_L$, and $C_{1,2}$ sufficiently small such that
for all $F \in \cU$ the corresponding lift $\tF$ preserves $C$-Lipschitz graphs with
constants $\ell_0, \la_0, C_{2,1}, C_{1,2}$ and the following inequalities also hold:
$$ C > (\ell - \la)^{-1} \| \tau \|_L \geq \| \rho^{\pm}_{\Alt} \|_L,$$
$$\widehat{\la} = \frac{\la_0 \ell_0 + C_{2,1} C_{1,2}}{\ell_0 - C_{1,2} C} < 1.$$

Let $\trho^+$ be the lift of $\rho^{+}_{\Alt}$.  By Lemma~\ref{tfie-l},
\begin{eqnarray*}
\| (\cT^+_F)^n (\trho^+) - \trho^+ \|_\infty & \leq &  \sum^{n-1}_{k=0} \|  (\cT^+_F)^{k+1} (\trho^+) - (\cT^+_F)^{k} (\trho^+)
\|_\infty \\
& \leq &  (\sum^{n-1}_{k=0 }\widehat{\la}^k) (1+ C) \delta \leq \frac{1+C}{1-\widehat{\la}} \delta.
\end{eqnarray*}
Choosing $\delta >0$ so that $\frac{1+C}{1-\widehat{\la}} \delta = \epsilon$ it follows that
$\trho^+_F = \lim (\cT^+_F)^n (\trho^+)$ is $\epsilon$-close to $\trho^+$.
Similarly, we obtain that $\trho^-_F$ is $\epsilon$-close to $\trho^-$.
\hfill
$\Box$

\section{Examples.}
\label{example-s}

\subsection{Perturbation of affine maps with annular attractors}
In the previous section we showed that the upper and lower boundaries
of the attractor $\Omegalt$ of an affine map $\Alt$ moves 
continuously under Lipschitz perturbations. 
In view of Lemma~\ref{annulus-characterization-l}, Remark~\ref{strict-r}
and Theorem~\ref{strong-affine-th}
we obtain the following result.

\begin{proposition}
\label{perturbation-p}
Let $\tau : \T \rightarrow \R$ be a Lipschitz function. 
Then there exists $\la_0 \in (1/\ell, 1)$ such that for any $\la \in [\la_0,1)$
and for all $F$ in an appropriate neighborhood $\cU$ of $\Alt$
in $\lip (\overline{U}, U)$ we have that  $\Omega_F$ is an annulus 
where $U = \T \times (-T_0, T_0)$ is such that $\Alt \in \lip (\overline{U}, U)$.
\end{proposition}

We will apply the above proposition to exhibit annular attractors in
an explicit family of endomorphisms of $\C^{*} = \C \setminus \{0\}$.
More precisely, we consider the family 
$$\begin{array}{rccc}
f_{\la,c} :& \Cstar & \rightarrow & \C \\
& z  & \mapsto & (\la |z| + 1- \la)\frac{z^2}{|z|^2} +c,
\end{array}$$
where $c \in \C$ and $\la \in (0,1)$. 
Observe that 
$f_{\la,0} (\Cstar) = \{ z \in \C \,\,/\, |z| > 1- \la \} \subset \Cstar$ and
$f_{\la,0}$ acts as multiplication by $2$ on the arguments and as the affine contraction
$\la | z| + 1- \la$ on radial lines. 
Also,  $f_{\la,0}(\S1) = \S1$ where $\S1=  \{ | z | =1 \}$.
The  map  $f_{\la,c}$ may be written as the postcomposition of $f_{\la,0}$ by the translation $z \mapsto z+c$.
Therefore, $f_{\la,c} (\Cstar) = \{ z \in \C \,\,/\, |z| > 1- \la \} +c$.

This family $f_{\la,c}$ is closely related to the well known and extensively studied
quadratic family $Q_c(z) = z^2 +c$.  The action of $Q_0$ as $|z|^2$ on radial
lines has been replaced by an affine contraction.

Our main interest here will be on maps $f_{\la,c}$ for $|c|$ small.
For $|c| < 1- \la$, by the above considerations we have that
$f_c (\Cstar) \subset \Cstar$. Moreover, it is easy to check that
$$V = \{ z \in \Cstar \,\,/\,    \exp(-\frac{|c|}{1-\la}) < |z| <   \exp(\frac{|c|}{1-\la} ) \} $$
is a trapping region for $f_{\la,c}$ when $|c|$ is sufficiently small. That is, $f_{\la,c} (\overline{V}) \subset V$.
We say that $$\Omega_{\la,c} = \cap_{n \geq 0} f^n_{\la, c} (\overline{V})$$
is the attractor for the dynamics of $f_{\la, c}$. 

\begin{proposition}
\label{family-e}
There exists $\la_0 \in (0,1)$ such that if $ \la \in [\la_0,  1 )$,
then  the attractor of $f_{\la ,c}$ is
homeomorphic to a closed topological annulus
for all $c$ in a punctured neighborhood of the origin (which depends on $\la$).
\end{proposition}

The proof of this fact relies on considering an appropriate rescaling
of $f_{\la, c}$ as $c$ goes to $0$. Namely, consider 
$$\begin{array}{rccc}
\iota :& \cyl & \rightarrow & \Cstar \\
& \tht & \mapsto &  \exp ( 2 \pi (t + i \theta))
\end{array}$$
and for $|c| < 1 - \la$, let
$$\widehat{f}_{\la, c} = \iota^{-1} \circ f_{\la, c} \circ \iota.$$

\begin{lemma}
\label{reescaling-l}
For $\eta >0$ let $h_\eta \tht = (\theta, \eta t)$.
Then, for all $\alpha \in \T$,  $$h^{-1}_\eta \circ \widehat{f}_{\la, \eta exp(2 \pi i \alpha)} \circ h_\eta$$
converges, as $\eta \rightarrow 0$,  in the $\cC^1$ topology to:
$$F_{\la,\alpha} \tht = ( 2 \theta, \la t + \frac{1}{2 \pi} \cos 2 \pi (\alpha - 2\theta)).$$ 
\end{lemma}

\proof
Fix $\la \in (0,1)$.
Let $c(\eta) = \eta \exp(2 \pi i \alpha)$ and
$$(\Phi_\eta \tht, \Psi_\eta \tht) =  \widehat{f}_{\la, c(\eta)} \circ h_\eta \tht.$$
Therefore,
$$h^{-1}_\eta \circ \widehat{f}_{\la , \eta \exp(2 \pi i \alpha)} \circ h_\eta
= (\Phi_\eta, \frac{\Psi_\eta}{\eta}).$$
We must show that $\Phi_\eta \tht \rightarrow 2 \theta$ and
that $$\frac{\Psi_\eta \tht}{\eta} \rightarrow \la t + (2 \pi)^{-1} \cos (\alpha - 2 \theta)$$ in the $\cC^1$ topology, as $\eta \rightarrow 0$.
>From $f_{\la, c} \circ \iota \circ h_\eta = \iota \circ (\Phi_\eta, \frac{\Psi_\eta}{\eta})$ we obtain:
\begin{equation}
\label{family1-e}
(\la \exp(2 \pi \eta t) + 1- \la) \exp(2 \pi i 2 \theta) + \eta \exp(2 \pi i \alpha)
\hfill  = \exp (2 \pi (\Psi_\eta \tht + i \Phi_\eta \tht)).
\end{equation}
It follows that
\begin{equation}
\label{family2-e}
\exp( 2 \cdot 2 \pi i (\Phi_\eta \tht - 2 \theta)) =
\frac{\la \exp(2 \pi \eta t) + 1 -\la + \eta \exp(2 \pi i (\alpha - 2 \theta))}{\la \exp(2 \pi \eta t) + 1 -\la + \eta \exp(-2 \pi i (\alpha - 2 \theta))}
\end{equation}
converges to the constant function $1$ in the $\cC^1$ topology
for maps from $\cyl$ into $\R$.
Therefore,
$$\Phi_\eta \tht - 2 \theta \rightarrow 0$$
in the $\cC^1$ topology for maps from $\cyl$ into $\T$.
Now
\begin{eqnarray}
\label{family3-e}
& & \hspace{-2cm} 
\exp( 2 \cdot 2 \pi i (\Phi_\eta \tht - 2 \theta))^{1/\eta} \nonumber \\
&=& \left(\frac{1 + \eta(2 \pi \la t  +  \exp(2 \pi i (\alpha - 2 \theta)) + O(\eta)}{1 + \eta(2 \pi \la t  +  \exp(-2 \pi i (\alpha - 2 \theta)) + O(\eta)}\right)^{1/\eta} \nonumber \\
& \rightarrow & \frac{\exp(2 \pi \la t) + \exp(2 \pi i ( \alpha - 2 \theta)}{\exp(2 \pi \la t) + \exp(-2 \pi i ( \alpha - 2 \theta)} \mbox{ as } \eta \rightarrow 0 \nonumber \\
& = & \exp ( 2 \cdot i \sin 2 \pi (\alpha - 2 \theta))
\end{eqnarray}
where the limit is in the $\cC^1$ topology.
>From (\ref{family1-e}),
$$\exp(2 \pi \frac{\Psi_\eta \tht}{\eta}) 
= (1 + \eta(2 \pi \la t + \exp(2 \pi i (\alpha - 2 \theta) + O(\eta))))^{1/\eta} \cdot \exp(2 \pi i (2 \theta - \Phi_\eta \tht))^{1/\eta},$$
which, in view of (\ref{family3-e}) converges to
$$\exp(2 \pi \la t + \exp(2 \pi i (\alpha - 2 \theta))) \cdot \exp(- i \sin 2 \pi (\alpha- 2 \theta)) = \exp(2 \pi \la t+ \cos 2 \pi (\alpha - 2 \theta))$$
in the $\cC^1$ topology.
It follows that $$ \frac{\Psi_\eta}{\eta} \rightarrow 
 \la t+ \frac{1}{2 \pi} \cos 2 \pi (\alpha - 2 \theta)$$
which establishes the Lemma.
\hfill
$\Box$

\medskip
>From Proposition~\ref{perturbation-p}, there exists $\la_0$ such that for all $\la \in [\la_0,1)$ there exists a neighborhood
$\cU$ of $\{ F_{\la, \alpha} \,\,/ \, \alpha \in \T \} \subset \lip(\overline{U}, U)$ so that the attractor $\Omega_F$ is 
an annulus for all $F \in \cU$ where $U= \T \times (-(1-\la)^{-1}, (1-\la)^{-1})$. Since  for $\eta$ sufficiently small,
say $0< \eta < \eta_0$, and  for all $\alpha \in \T$,
$$h^{-1}_\eta \circ \widehat{f}_{\la, \eta exp(2 \pi i \alpha)} \circ h_\eta \in \cU,$$
  it follows that the attractor set of 
$\widehat{f}_{\la, \eta exp(2 \pi i \alpha)} \in \lip(\overline{W},W)$ is an annulus, where $W= \T \times (-(1-\la)^{-1} \eta,
(1-\la)^{-1} \eta)$. Therefore, the attractor of $f_{\la, \eta exp(2 \pi i \alpha)}$ is an annulus for all $0 < \eta < \eta_0$ and
all $\alpha \in \T$ which proves  the claim of Proposition~\ref{family-e}.

\subsection{With interior and not an annulus.}

\begin{proposition}
\label{fat-hole-ex}
For all $\la >1/2$ there exists a $\cC^\omega$ function $\tau: \T \rightarrow \R$ such that
the attractor set $\Omegalt$ of $\Alt \tht = (2 \theta, \la t + \tau (\theta))$ has
non-empty interior and it is not a topological annulus. Moreover, this property is 
robust. That is, let $U = \T \times (-T_0,T_0) \subset \cyl$ be 
such that $\Alt(\overline{U}) \subset U$.
Then there exists a neighborhood $\cU$ of $\Alt$ in $\lip(\overline{U},U)$ such that for all
$F \in \cU$ the attractor set $\Omega_F$  has
non-empty interior and it is not a topological annulus.
\end{proposition}

A similar example can be constructed for all $\la > 1/\ell$ and $\ell >2$. Here we specialize in the case $\ell =2$ for 
the sake of simplicity of the exposition. 

\medskip
\noindent
{\bf Construction of the example.}
We endow $\T$ with its standard orientation  and use interval notation accordingly.
Let $p \geq 1$ be such that $\la + \cdots + \la^{p-1} > 1$ and consider $\eta < 1$ such that:
\begin{equation}
\label{ex1-e}
(\la + \cdots + \la^{p-1}) \eta > 1.
\end{equation}
Consider a period $p$ periodic cycle $ \theta_0, \dots, \theta_{p-1} = \m^{p-1}_2 (\theta_0)$ with subindices 
$\operatorname{mod} p$ and respecting cyclic order such that $|I_0| = |(\theta_0,\theta_1)| > 1/2$
where $|I|$ denotes the length of the interval $I \subset \T$.  For example, let $\theta_1 = 1/(2^p-1)$.

Consider positive constants $T_0, T_1,  \delta, \epsilon_0, \dots, \epsilon_{p-1}, \la^\prime$ and 
a natural number $N \geq 2$ such that
\begin{eqnarray}
T_0 & > & \frac{1}{1-\la}, \label{ex2-e} \\
\frac{T_0}{\la} & > & \frac{\la - \la^N}{1-\la} + \frac{\la^N T_0}{1-\la} = T_1,  \label{ex3-e}\\
\delta & = & | I_0| - 1/2, \label{ex4-e} \\
0 & < & \epsilon_0 < \delta \cdot 2^{-N+1}, \label{ex5-e} \\
2 \epsilon_j & < & \epsilon_{j+1} \mbox{ for } j= 1, \dots,p-1, \label{ex6-e} \\
\la \eta & < & \la^\prime < \la. \label{ex7-e}
\end{eqnarray}
Let $\tau : \T \rightarrow \R$ be a $\cC^\infty$ function such that:
\begin{eqnarray}
\tau(\theta) & =  & 0 \mbox{ for all } \theta \in I^{\epsilon_0}_0 \label{ex8-e} \\
\tau(\theta) & = & \la^\prime  \mbox{ for all } \theta \in I^{\epsilon_j}_j    \label{ex9-e} 
\end{eqnarray}
where $I^{\epsilon_j}_j = (\theta_j + \epsilon_j , \theta_{j+1}-\epsilon_{j})$,
\begin{eqnarray}
\tau(\theta_j) &  = & T_0 \mbox{ for } j= 1, \dots,p-1, \label{ex10-e} \\
0  <  \tau(\theta) & \leq & T_0  \mbox{ for } \theta \in \T \setminus I_0. \label{ex11-e}
\end{eqnarray}

The next three lemmas are devoted to prove
that  the properties stated in Proposition~\ref{fat-hole-ex} hold for a function $\tau: \T \rightarrow \R$
as above. It follows that there is a Lipschitz close $\cC^\omega$ function for which the statement of 
the proposition holds.

\begin{lemma}
\label{first-ex-l}
Let $\rho^\pm = \rho^\pm_{\Alt}$ on $\T$. Then

{\em (1)} $\rho^- (\theta)= 0 $ for all  $\theta \in \T$.

{\em (2)} $\rho^+ (\theta_j) = (1- \la)^{-1} T_0$ for all $j=0, \dots, p-1.$

{\em (3)} $\rho^+ (\theta_0 + \frac{1}{2}) < \la^{-1} T_0$. 
\end{lemma}

\proof
(1) Since $$|I^{\epsilon_0}_0| = |I_0| - 2 \epsilon_0 > |I_0| - 2 \delta \cdot 2^{-N+1} \geq |I_0| - \delta,$$
it follows that $\m_2 (I^{\epsilon_0}_0) = \T$. Therefore, given $\thetap_0 \in \T$ and $n \geq 1$ there exists 
$\thetap_{n} \in I^{\epsilon_0}_0$ such that $2 \thetap_{n} = \thetap_{n-1}$. Hence, $\rho^- (\thetap_0) \leq 
t_\la ( (\thetap_n)) = 0$ (see Proposition~\ref{affine-semiconjugacy-p}).  Also $ 0 \leq  \rho^- (\thetap_0)$ 
since $\tau \geq 0$.

\smallskip
(2) Since $\tau \leq T_0$, it follows that $\rho^+ \leq (1-\la)^{-1} T_0$.
Now  $\rho^+ (\theta_j) = (1-\la)^{-1} T_0$ because  each one of the periodic points $\theta_j$ has as an infinite
backward orbit along the
periodic orbit $\theta_0, \dots, \theta_{p-1}$.

\smallskip
(3) Let $\thetap_0 = \theta_0 + 1/2$ and $(\thetap_n) \in \cS$ be a backward orbit. Let 
$$N_0 = \min \{ n \,\,/\, \thetap_n \notin I^{\epsilon_0}_0 \cup \cdots \cup I^{\epsilon_{p-1}}_{p-1}\}.$$
Then $$\Tdist (\thetap_{N_0}, \{ \theta_j \}^{p-1}_0)  
\leq  \max \{ \epsilon_j \} = \epsilon_0 < \delta \cdot 2^{-N+1}.$$
Hence, $$\delta =  \Tdist (\thetap_{0}, \{ \theta_j \}^{p-1}_0) <  \delta \cdot 2^{-N+1} \cdot 2^{N_0}.$$
We conclude that $N_0 \geq N$ and
\begin{eqnarray*}
\rho^+(\thetap_0) & \leq& \tau(\thetap_1) + \la \tau(\thetap_2) + \cdots +  \la^{N_{0}-2} \tau(\thetap_{N_0-1})
+ \frac{\la^{N_0} T_0}{1-\la} \\
                  & \leq & \la + \la^2 + \cdots + \la^{N_{0}-1} +  \frac{\la^{N_0} T_0}{1-\la} \\
                  & =  & \frac{\la - \la^{N_0}}{1-\la} + \frac{\la^{N_0} T_0}{1-\la} \\
                  & \leq & \frac{\la - \la^{N}}{1-\la} + \frac{\la^{N} T_0}{1-\la}  < \la^{-1} T_0.
\end{eqnarray*}
\hfill $\Box$

\begin{lemma}
\label{non-annular-l}
There exists an open neighborhood $\cU$ of $\Alt$ in  $\lip(\overline{U}, U)$ such
that for all $F \in \cU$ the attractor $\Omega_F$ is not a closed topological annulus. 
\end{lemma}

\proof
If we denote by $\gamma^\pm$ the graphs of $\rho^\pm_{\Alt}$, then
\begin{eqnarray*}
& & \hspace{-2cm} \min \{ t \,\, /\, (\theta_1, t) \in \Alt (\gamma^+) \} \\
 & = & \min \{ \la \rho^+ (\theta_0 + 1/2) + 
\tau (\theta_0 + 1/2), \la \rho^+ (\theta_0) + \tau (\theta_0) \} \\
&=& \min \{ \la \rho^+ (\theta_0 + 1/2) , (1-\la)^{-1} T_0 \} \\
& < & T_0 \\
&=& \max \{ 0, T_0 \} \\
&=& \max \{ \la \rho^- (\theta_0 + 1/2) + \tau (\theta_0 + 1/2), \la \rho^- (\theta_0) + \tau(\theta_0) \} \\
&=& \max \{ t \,\, /\, (\theta_1, t) \in \Alt (\gamma^-) \}.
\end{eqnarray*}
By Lemma~\ref{annularF-l}, it follows that for all $F$ close to $\Alt$ the attractor $\Omega_F$ is not 
an annulus.
\hfill
$\Box$

\begin{lemma}
\label{interior-l}
There exists an open neighborhood $\cU$ of $\Alt$ in  $\lip(\overline{U}, U)$ such
that for all $F \in \cU$ the attractor $\Omega_F$ has non-empty interior.
\end{lemma}

\proof
Let $\eta^\prime$ be such that
$$(\la + \cdots + \la^{p-1}) \eta > \eta^\prime >1.$$
By Proposition~\ref{uplow-p}
we may assume that for all $F$ sufficiently close to
$\Alt$ the lower boundary $\rho^-_F : \T \rightarrow \R$ is well defined
and 
$$\rho^-_F <  (\la^\prime)^2.$$
Consider the subsets $R^F_1$ and $R_2$ of $\cyl$ defined by:
$$R^F_1 = \{ \tht \,\,/ \, \rho^-_F(\theta) \leq t \leq \la \},$$
$$R_2   = \overline{I^{\epsilon_0}_0} \times [\la, \eta^\prime] \cup \cup^{p-1}_{j=2 } 
(\overline{I^{\epsilon_j}_j} \times [ \la, (\la + \cdots + \la^j) \eta]).$$
Let $R^F = R^F_1 \cup R_2$. We will show that the $F(R^F) \supset R^F$ for all $F$ in a sufficiently small
neighborhood of $\Alt$. 

To simplify notation, let $R = R^{\Alt}$.
We start by showing that $\Alt ( \operatorname{int} R ) \supset R_2$. In fact, since for $j=1, \dots, p-1$ 
$$ I^{\epsilon_j}_j \times ( 0, (\la + \cdots + \la^j) \eta) \subset \operatorname{int} R$$
and $\tau( I^{\epsilon_j}_j) =\{ \la^\prime\}$,  it follows that 
$$ \m_2 (I^{\epsilon_j}_j) \times (\la^\prime + \la(\la + \cdots + \la^j) \eta) \subset
\Alt(\operatorname{int} R).$$
By (\ref{ex6-e}), $\m_2 (I^{\epsilon_j}_j) \supset  \overline{I^{\epsilon_{j+1}}_{j+1}}$ and,
by (\ref{ex7-e}), $ (\la^\prime + \la(\la + \cdots \la^j) \eta) \supset 
[\la, \la \eta + \la (\la + \cdots + \la^j) \eta]$.
Therefore,  $\Alt ( \operatorname{int} R ) \supset R_2$.
Hence, for all $F$ sufficiently close to $\Alt$ we also have
that 
\begin{equation}
F(\operatorname{int} R^F ) \supset R_2 \label{contains-r2-e}
\end{equation}
since the boundaries of $R^F$ move continuously with $F$.

For $F$ sufficiently close to $\Alt$ the circle $\T \times \{  \la \}$ has image 
contained in $t > \la^\prime \la$. Therefore, 
\begin{equation}
F(\{\tht \,\, / \, \rho^-_F(\theta) \leq t \leq \la \} \supset 
\{\tht \,\, / \, \rho^-_F(\theta) \leq t \leq \la \la^\prime \}. \label{contains-bottom-e}
\end{equation}

Since $\Alt ( I^{\epsilon_0}_0 \times (\la \eta, \eta^\prime)) \supset
\T \times (\la \la \eta, \la \eta^\prime) \supset \T \times [\la \la^\prime, \la]$,
for all $F$ sufficiently close to $\Alt$ we also  have that 
\begin{equation}
F(R^F) \supset \T \times [\la \la^\prime, \la]. \label{contains-top-e}
\end{equation}
>From (\ref{contains-bottom-e}) and (\ref{contains-top-e}) we conclude that $F(R^F) \supset R^F_1$.
By (\ref{contains-r2-e}), we obtain that $F(R^F) \supset R^F$.

Since for all $F$ sufficiently close to $\Alt$ we have that $F(R^F) \supset R^F$, it follows 
that $R^F \subset \Omega_F$ and therefore $\Omega_F$ has non-empty interior.
\hfill
$\Box$

\bibliographystyle{plain}

\end{document}